\documentclass[11pt]{amsart}

\usepackage[mathscr]{eucal}
\usepackage{ulem}
\usepackage{stmaryrd}
\usepackage{amssymb}
\usepackage{amsfonts}
\usepackage{amsmath}
\usepackage{amsthm}

\usepackage{graphicx,color}

\usepackage[all,color,frame,import]{xy}
\usepackage{epsfig}
\usepackage{warmread}

\newtheorem{theorem}{Theorem}[section]
\newtheorem{lemma}[theorem]{Lemma}

\newtheorem{definition}[theorem]{Definition\rm}

\newcommand{\git}{/\! \!/}
\newcommand{\hk}{/\! \!/\!\!/}
\newcommand{\CP}{\mathbb{CP}}
\newcommand{\proj}{\mathbb{P}}
\newcommand{\C}{\mathbb{C}}
\newcommand{\R}{\mathbb{R}}

\newcommand{\HK}{hyper-K\"ahler }
\newcommand{\K}{K\"ahler }
\newcommand{\nrm}[1]{|\!| #1|\!|}
\newcommand{\abs}[1]{|#1|}

\newcommand{\res}{\text{Res}}
\newcommand{\JKres}{\text{JKRes}}
\newcommand{\jkres}{\text{jkres}}
\newcommand{\mres}{\text{res}}
\newcommand{\gies}[2]{\ensuremath{{\mathcal{M}}^{g}_{{#1} ,{#2} }}} %moduli spaces of torsion free sheaves, gieseker space
\newcommand{\uhl}[2]{\ensuremath{{\mathcal{M}}^u_{#1,#2}}} %moduli spaces of (ideal) instantons, uhlenbeck space
\newcommand{\modu}[2]{\ensuremath{{\mathcal{M}}^o_{#1,#2}}} 
\newcommand{\adhm}[2]{\ensuremath{{\mathbb{M}}_{{#1} ,{#2} }}}

\newcommand{\csympladhm}[2]{\ensuremath{\left(\mathbb{M}^{Sp}_{#1,#2}\right)_{\lambda}}}
\newcommand{\sympladhm}[2]{\ensuremath{{\mathbb{M}}^{Sp}_{{#1},{#2} }}}
\newcommand{\sympl}[2]{\ensuremath{{\mathcal{M}}^{Sp}_{{#1},{#2} }}}
\newcommand{\sportadhm}[2]{\ensuremath{{ \mathbb{M}}^{SO}_{{#1},{#2} }}}
\newcommand{\sport}[2]{\ensuremath{{ \mathcal{M}}^{SO}_{{#1},{#2} }}}
\newcommand{\im}{\text{im }}
\renewcommand{\emph}{\normalem} %%makes \emph behave normally when using the ulem package
\newcommand{\ch}{\ensuremath{\text{ch}}}
\newcommand{\volu}{\ensuremath{\text{vol}}}
\renewcommand{\coprod}{\sqcup}

\title[Equivariant Volumes of  quotients \& Instanton Counting]{Equivariant volumes of non-compact quotients and instanton counting}
\author{Johan Martens}
\address{Department of Mathematics, University of Toronto, Toronto Ontario M5S 2E4, Canada}%\curraddr{\textsc{Max-Planck-Institut f\"ur Mathematik, Vivatsgasse 7, 53111 Bonn, Germany}}
\email{jmartens@math.toronto.edu}\date{\today}
\begin{document}

\begin{abstract}
Motivated by Nekrasov's instanton counting, we discuss a method for calculating equivariant volumes of non-compact
quotients
in symplectic and \HK geometry by means of the Jeffrey-Kirwan
residue formula of non-abelian localization.  In order to
overcome the non-compactness, we use varying symplectic cuts to reduce the problem to a compact setting, and study
what happens in the limit that recovers the original problem.  We implement this method for the ADHM construction of the moduli spaces of framed
Yang-Mills instantons on $\R^{4}$ and rederive the formulas for the equivariant volumes obtained earlier by Nekrasov-Shadchin, expressing these volumes as iterated residues of a single rational function.
\end{abstract}

\maketitle

\section{Introduction}
In this paper we develop a method for calculating equivariant volumes of non-compact symplectic spaces (where they are also known as \emph{regularized volumes}) equipped with a Hamiltonian action of a torus $T$, in the case when these spaces are finite-dimensional symplectic reductions or \HK quotients.  The Hamiltonian action gives a moment map $\mu:M\rightarrow \mathfrak{t}^*$, and we are interested in the improper integrals  
\begin{equation}\label{hintmichael}\int_M e^{\langle \mu,t\rangle }\frac{\omega^n}{n!}\end{equation}
of the function $e^{\langle \mu,t\rangle }$ with respect to the
Liouville measure $\frac{\omega^n}{n!}$.
Under certain conditions this integral converges for $t$ in an
open subset of the Lie algebra $\mathfrak{t}_{\C}$.  Formally we may write it as 
$$\int_M e^{\omega + \langle\mu,t\rangle } $$ where $e^{\omega +\langle\mu,t\rangle}$  is a class in a  completion of the equivariant cohomology ring $H^*_T M$. 
Because $M$ is not compact, however, the integral of this class
cannot be directly defined cohomologically. 

The method described here is essentially a combination of the
symplectic cut of Lerman~\cite{cuts} and the Jeffrey-Kirwan
residue formula of non-abelian localization~\cite{nonabloc}.
The symplectic cut allows us to reduce the integral to an
integral over a compact space $M_{\lambda}$, such that the original integral is recovered in the limit as $\lambda$ goes to infinity.  Once in the compact setting, the residue formulas (which require compactness) allow the integral over the quotient to be expressed in terms of equivariant data on the space that one takes the quotient of.  The key observation is that the cutting operation introduces new fixed points `at infinity', whose contribution is crucial. 

In the particular case of symplectic reductions or \HK
quotients of vector spaces by linear actions, we are able to
describe the outcome of this procedure in greater detail, and
express the corresponding volumes as iterated residues of a
single rational function.  For \HK quotients for instance we have $$\int_{V\hk G} e^{\omega+\mu}=\frac{\tilde{K}}{\abs{W_G}}\res_+^{X_1,\ldots, X_n} \frac{\varpi^2 e(\mu_{\C})}{\prod_j \rho_j}.$$ In the rational function $\frac{\varpi^2 e(\mu_{\C})}{\prod_j \rho_j}$ in question, $\varpi$ is the product of the positive roots of $G$, and both $e(\mu_{\C})$ and the denominator are products of linear factors entirely determined by the action of $G$ on $V$.  The iterated residue $\res_+^{X_1,\ldots, X_n} $ extracts the relevant information out of the rational function; its definition is given in section~\ref{linact}.

Our motivating examples come from physics: the equivariant
volumes of moduli spaces of framed instantons, which are the main ingredient in the recent instanton counting of Nekrasov~\cite{nekra}.  The famous ADHM construction~\cite{ADHM} realizes these spaces as \HK quotients, and using this we completely carry out the calculational strategy described below for them. 

We would like to emphasize that we interpret the integrals we study as \emph{regularized volumes} $\int e^{\omega+\mu}$, i.e. bona fide improper integrals, firmly rooted in symplectic geometry, which is true to the initial formulation of the problem of instanton counting in physics.  This is somewhat in contrast to other authors~\cite{naka,rndprt}, who work with a formal integral $\int 1$ in equivariant cohomology of the class 1, essentially defined through a localization formula.  From our more naive point of view, all localization formulas are theorems rather than definitions, and more importantly for us, this perspective allows for our calculational method, as the moment map gives the exact recipe for handling the residues of the new fixed points introduced by the cut.  

\subsection{Organization} The rest of this paper is organized
as follows.  In section~\ref{pw} we discuss the basic object of our study, the regularized volumes of non-compact symplectic spaces with Hamiltonian torus action, and review the Prato-Wu theorem that calculates their volume.  In section~\ref{basisidee} we discuss the Jeffrey-Kirwan residue theorem, and develop our main conceptual idea, combining the residue theorem with symplectic cutting to calculate the equivariant volume of quotients.  In section~\ref{linact} we specialize to \emph{linear} symplectic reductions and \HK quotients, and discuss the specific implementation of our calculational method in these cases.  
As examples of this calculational method we implement it in section~\ref{nstntn} for the ADHM spaces for the classical gauge groups, and compare the results with the work of Nekrasov-Shadchin~\cite{ABCD} whose formulas we rederive. 
\subsection{Notation} We will use the following notation throughout:
$$\begin{array}{cl}
\langle .,.\rangle & \text{natural pairing between
}\mathfrak{t}^*\text{ and }\mathfrak{t}\\ &\text{or their
  complexifications }\mathfrak{t}^*_{\C}\text{ and
}\mathfrak{t}_{\C} \\ \ & \ \\ M\git_{\lambda}G &
\text{symplectic reduction of }M\text{ by }G\\ \ &\text{at a
  central value }\lambda\text{ of the moment map} \\ \ & \ \\
M_{\lambda}, V_{\lambda} &\text{symplectic cut of }M,V\text{ with respect to a circle action}\\ &\text{at value }\lambda\in \mathfrak{u}(1)^*\\ \ & \ \\ T^n \ni (e_1,\dots,e_n) & \text{elements of a split torus $T=U(1)^n$, and} \\ = (e^{\epsilon_1},\ldots,e^{\epsilon_n})  & \text{corresponding coordinates on Lie algebra $\mathfrak{t}$;}\\ & \text{we will use corresponding Greek and Roman} \\ & \text{lower-case letters for elements of }\mathfrak{t}\text{ and }T\text{ respectively}
\end{array}$$
\subsection*{Acknowledgements.}
I would like to thank Lisa Jeffrey, Mikhail Kogan, Nikita Nekrasov, Mich\`ele Vergne and Peter Woit for useful conversations.  The bulk of the work presented here was done while I was a graduate student at Columbia University, and I would like to thank my advisor Michael Thaddeus for continuous guidance and encouragement, as well as a careful reading of the final version of this article.  Finally I would like to add that this work was largely inspired by reading the early paper~\cite{higgs-int}.  This article was partially written during a stay at the Max-Planck-Institut f\"ur Mathematik, Bonn, whose hospitality and support is gratefully acknowledged.

\section{The Duistermaat-Heckman theorem for non-compact spaces}\label{pw}
The famous formula of Duistermaat-Heckman~\cite{DH2} establishes an expression for the integral of the
\emph{equivariant volume} of a compact symplectic manifold $(M,\omega)$ equipped with a Hamiltonian action of a torus
$T$ with moment map $\mu$:
 \begin{equation*}\int_{M}e^{\omega+\mu}=\sum_{F\subset M^{T}} \int_{F } \frac{e^{\omega+\mu}}{e_{T}(\nu_{F})}\end{equation*}
 where the sum in the right hand side is over the connected components of the fixed point set for the action of $T$, $\nu_{F}$ is the
normal bundle of $F$ in $M$, and $e_{T}(\nu_{F})$ is its equivariant Euler class.  In~\cite{atibott,bv} this formula was 
established as a particular instance of localization in
equivariant cohomology, using the Cartan model. In~\cite{prato-dh},
this theorem is
discussed in the case where $M$ is no longer compact:
\begin{theorem}[Prato-Wu] Let $(M,\omega)$ be a (non-compact) symplectic manifold equipped with a Hamiltonian action of a torus $T$. If the fixed-point set $M^T$ is compact and there exists a $t_{0}\in\mathfrak{t}=\text{Lie}(T)$ such that the corresponding component of the moment
map $\mu_{t_{0}}=\langle\mu,t_0\rangle$
is proper and bounded on one side (assume bounded above) then the Duistermaat-Heckman formula \begin{equation}\label{pwform}
\int_M e^{\langle \mu,t\rangle }\frac{\omega^n}{n!}=\sum_{F\subset M^{T}} \int_{F } \frac{e^{\omega+\langle\mu,t\rangle}}{e_{T}(\nu_{F})}
\end{equation}
still holds true, provided one interprets the right-hand side as the integral of the function $e^{\langle\mu, t \rangle}$ with
respect to the measure induced by the Liouville volume $\frac{\omega^{n}}{n!}$, and one restricts oneself to $t$
inside of an certain open cone $C\subset\mathfrak{t}_{\C}$.
\end{theorem}
The open cone occurring in the theorem is exactly the set of all $t$ for which the integral in the left hand side of (\ref{pwform}) converges.  The right-hand side can now be interpreted as a function on $\mathfrak{t}_{\C}$ with poles on an
arrangement of hyperplanes containing the origin, and the cone used is contained in a component of the complement of
the hyperplanes.  

One can prove this theorem by means of the \emph{symplectic cut} construction introduced by Lerman~\cite{cuts}: the
conditions stated above guarantee the existence of a circle $U(1)$ in $T$ whose moment map $\mu_{c} $ is proper and bounded above, hence the symplectic cut $$
M_{\lambda} = \left(M\times \C\right)\git_{\lambda} U(1) = \mu_{c}^{-1}(-\infty,\lambda) \coprod M\git_{\lambda} U(1)$$ with respect to this circle
is a compact symplectic space, still equipped with a $T$-action.  As the open dense subset inherits its symplectic structure from the inclusion of 
$\mu_{c}^{-1}(-\infty,\lambda)$ in $M$, one can moreover recover the integral over $M$ we are interested in as the limit
of integrals over the symplectic cuts:
\begin{equation}\label{blah}\int_{M}e^{\omega+\mu}=\lim_{\lambda\rightarrow \infty} \int_{M_{\lambda}}e^{\omega+\mu}. \end{equation}
For the compact spaces $M_{\lambda}$ the original Duistermaat-Heckman formula holds, and as $M^{T}$ is assumed to be
compact, one can remark that the fixed point set $\left(M_{\lambda}\right)^{T}$, for $\lambda$ large enough, consists
of the `original' fixed point components of $M$, together with `new' ones introduced by the cut:
$$
M_{\lambda}^T = \left(\mu_{c}^{-1}(-\infty,\lambda)\right)^T \coprod \left(M\git_{\lambda} U(1)\right)^T = M^T \coprod \left(M\git_{\lambda} U(1) \right)^T$$  
As the components of $M^T$ don't depend on $\lambda$, (\ref{blah}) becomes \begin{equation}\label{cmp}\int_M e^{\omega+\mu} = \sum_{F\subset M^{T}} \int_{F } \frac{e^{\omega+\mu}}{e_{T}(\nu_{F})}\ + \lim_{\lambda\rightarrow\infty} \ \sum_{F'\subset M\git U(1)^{T}} \int_{F'} \frac{e^{\omega+\mu}}{e_{T}(\nu_{F'})}. \end{equation}
Moreover,  in the
limit as $\lambda\to\infty$, one can easily show that the contributions of the new fixed points (i.e. the $T$-fixed points of $M\git_{\lambda} U(1)$) will vanish, $$\lim_{\lambda \rightarrow \infty} \ \sum_{F'\subset M\git U(1)^{T}} \int_{F'} \frac{e^{\omega+\langle\mu, t\rangle }}{e_{T}(\nu_{F'})} = 0,$$
provided that one exactly restricts to values of $t$ in the open cone $C$, which establishes (\ref{pwform}).  Indeed, the values of any new fixed-point component $F$ under the moment map $\mu$ will vary with the parameter $\lambda$, and the cone $C$ is such that the exponent of $e^{\langle\mu(F),t\rangle}$ is negative. It was pointed out to the author
by M. Vergne that
the original proof of \cite{prato-dh} essentially follows this line, although the technology of symplectic cuts had not
been introduced at the time.  

In \cite{paradan}, Paradan discusses the above formula in the much broader setting of
equivariant
cohomology with generalized coefficients.  For our purposes the naive interpretation given above suffices however,
and we will restrict ourselves to this framework.

\section{Equivariant volumes of non-compact quotients} \label{basisidee}
In this section we describe the main idea of this paper: in order to calculate the equivariant volume of a non-compact quotient $M\git G$, compactify both $M$ and $M\git G$ by a symplectic cut, use a residue formula to calculate the volume of the cut of $M\git G$, and recover the result in the limit.

\subsection{Residue Theorem}
Kirwan showed in~\cite{thesisfrances} that for the symplectic reduction of a compact symplectic space $(M,\omega)$ by a 
compact connected Lie group $G $ there is a surjective map, the \emph{Kirwan map}, $$\kappa:H^{*}_{G}(M,\mathbb{Q})\rightarrow H^{*}(M\git G,\mathbb{Q}).$$
In~\cite{nonabloc} Jeffrey and Kirwan established a method for calculating integrals over $M\git G$ of images of this
map, by means of the residue formula\footnote{Various appearances of this formula in the literature~\cite{jefkirriem,nonabloc2,arbitrary} differ in the constants used because of alternative notational conventions, e.g.\  the use of $\omega + i\mu$ instead of  $\omega+\mu$ as the equivariant symplectic form.  We will follow the version given in~\cite{arbitrary}. Other approaches to the residue formula where given in~\cite{notevergne,guilkalk}.}:
\begin{theorem}\label{lisafrances}
\begin{equation}\label{jefkir}
\int_{M\git G} \kappa (\alpha) e^{\omega} = \sum_{F\subset M^{T_G}} \frac{K}{\abs{W_G}}\ \JKres^{\Lambda}\ \varpi^{2}\int_{F}\frac{i^{*}_{F}\alpha
e^{\omega+\mu_{T_G}}}{e_{T_{G}}(\nu_{F})}
\end{equation}
where $T_G$ is a maximal torus of $G$, $K$ is the constant
$\frac{(-1)^{n_{\scriptscriptstyle +}}}{\volu(T_G)}$, $n_+$ is
the number of positive roots, the volume $\volu(T_G)$ is taken
with respect to an invariant metric on $T_G$, $W_{G}$ is the
Weyl group of $G$, $\varpi$ is the product of the positive roots of $G$, $i_F$ is the inclusion $i_F: F\hookrightarrow M$, 
and $\JKres^{\Lambda}$ is the Jeffrey-Kirwan residue operation.\end{theorem}
The residue $\JKres^{\Lambda}$ is a linear operation defined on a class of meromorphic forms
on a vector space, with poles on an arrangement of hyperplanes (i.e. the denominator has to split as a  product of linear factors).  The residue depends on the choice of a cone $\Lambda\subset
\mathfrak{t}$; for different choices of $\Lambda$ the residue will affect the terms in the right hand side
differently, but the total result is independent of the choice for $\Lambda$.  

The original definition of the residue was given in terms of
formal Fourier transforms, and a characterizing list of axioms was determined in~\cite[Proposition 8.11]{nonabloc}.  Another intrinsic formulation of the residue was given in~\cite{verbrion}.
We will however restrict ourselves to using another definition of the residue, given in~\cite{nonabloc2}, as iterated one dimensional residues.

Begin by looking at the one-dimensional case.  Let $f$ be a function on $\mathfrak{t}_{\C}\cong \C$ of the form $$f(x)=\sum_{j} g_j e^{\lambda_{j}x},$$ where $g_j$ are rational functions on $\C$, and $\lambda_j\in\R$.  Then define the residue of the $1$-form $f(x)dx$ to be \begin{equation}\label{bijna2}\jkres^+(f(x) dx)=\sum_{\lambda_j \geq 0} \sum_{b\in \C} \mres_b (g_j(x) e^{\lambda_j x}).\end{equation}
This one-dimensional residue will be the building block for defining higher-di\-men\-sio\-nal versions.  
A very useful little fact that we shall use often later on is:
if $p(x), q(x)$ are polynomials on $\mathfrak{t}_{\C}$ such that $\deg(p)+1\neq \deg(q)$, then \begin{equation}\label{bijna}\jkres^+ \frac{p(x)}{q(x)}dx=0, \end{equation} as can easily be seen by interpreting this rational form as a rational form on $\CP_1$, calculating $\jkres^+$ by using the Cauchy theorem, and looking at the residue at $\infty$.

Now the higher-dimensional version that we want to use is in practice obtained by iterating such one-dimensional residues $\jkres^+$.  Suppose $(X_1,\dots,X_n)$ is a linear coordinate system on $\mathfrak{t}$.  Then for any function given as a sum of terms of the form
\begin{equation} \label{goeiefunctie} f(z)=\frac{p(x) e^{\lambda(x)}}{\prod \rho_i(x) } \end{equation} with $q(x)$ a polynomial on $\mathfrak{t}_{\C}$, and $\lambda, \rho_i \in \mathfrak{t}^*$, we can still define one-dimensional residues $\jkres^+_{X_i}$ for any of the given coordinates, by treating all other coordinates as  generic constants.  The result of this is a new function of the form(\ref{goeiefunctie}), but this time only in the remaining variables $X_1,\ldots, \hat{X_i},\ldots, X_{n}$. 
We then follow~\cite{nonabloc2,locbycuts} in defining the residue operation $\JKres^{\Lambda}$ inductively as an iterated residue:
\begin{equation}\label{beupf}\begin{split}\JKres^{\Lambda} (h[dx]) =\triangle\jkres^+_{X_1}\ldots \jkres^+_{X_n}h(X_1,\ldots,X_n) [dX]_1^n \end{split}\end{equation}
for any generic coordinate system $(X_1,\ldots, X_n)$ such that $(0,\ldots,0,1)\in \Lambda$, where $\triangle$ is the determinant of any $n\times n$ matrix whose columns form an orthonormal basis of $\mathfrak{t}$ defining the same orientation as $(X_1,\ldots,X_n)$ (a metric on $\mathfrak{t}$ is tacitly understood).  Different choices for $\Lambda$ will give different results when applying  $\JKres^{\Lambda}$ to the summands in the RHS of (\ref{jefkir}), but the total result will be unaffected.  In general $\JKres^{\Lambda}$ will only give a non-zero contribution for the summands corresponding to components $F\subset M^{T}$ where $\mu(F)\in \Lambda^{*}$, the dual cone of $\Lambda$.  In~\cite{nonabloc,nonabloc2} it was shown that this definition of $\JKres^{\Lambda}$ coincides with the earlier one given in terms of a Fourier transform.

As indicated above the residue operation acts on meromorphic forms rather than functions, but for notational convenience we will suppress the relevant form part $[dX]_1^n$ throughout. 

\subsection{Equivariant volumes of symplectic reductions}
We want to apply theorem~\ref{lisafrances} for the calculation of equivariant volumes of non-compact spaces that arise as
symplectic reductions, in the sense described in section~\ref{pw}.  We will take the same approach as the one
outlined above for the theorem of Prato-Wu in (\ref{blah}): one can relate the integral on the non-compact space to a limit of
integrals on compact spaces by means of varying symplectic cuts, and then use the available residue theorems for the compact
spaces.  In sharp contrast however with the theorem of Prato-Wu is that the contributions of the `new' fixed
points introduced by the cuts will not vanish in the limit.  In fact, in the cases that we study later on, these
`fixed points at infinity' carry all the information. 

 In order to implement this we need the following straightforward
equivariant generalization of the Jeffrey-Kirwan theorem:
\begin{theorem}\label{equiJK}
Let $G$ and $H$ be two (compact, connected) Lie groups with commuting Hamiltonian actions on a compact symplectic manifold
$(M,\omega)$.  Then the action of $H$ descends to the symplectic reduction $M\git G$, and with the equivariant Kirwan map $\kappa:H^{*}_{G\times
H}M\rightarrow H^{*}_{H}(M\git G)$ and the notation as in theorem~\ref{lisafrances} we have
\begin{equation}
\int_{M\git G}\kappa(\alpha)e^{\omega +\mu_{H}}=\sum_{F\subset M^{T_{G}}} \frac{K}{\abs{W_G}}\JKres^{\Lambda}\varpi^{2}\int_{F}\frac{i^{*}_{F}\alpha
e^{\omega+\mu_{T_G}+\mu_H}}{e_{T_{G}}(\nu_{F})}.
\end{equation}
\end{theorem}

As our proof of this is a fairly straightforward reduction of the equivariant to the non-equivariant $H=1$ case and mainly consists of technical remarks, we postpone it to appendix~\ref{appen}.

Now suppose that we have two commuting Hamiltonian group actions on a symplectic manifold $(M,\omega)$, of a compact
connected Lie group $G$ and a torus $T$ with moment maps $\mu_{G}$ and $\mu_{T}$, and suppose that $\mu_{T}$ has a
component that is proper and bounded from below.  Then, as before, we can find a circle in $T$ such that the
symplectic cuts of both $M$ and $M\git G$ with respect to this circle are compact.  Moreover, since all group actions
commute and a symplectic cut is a symplectic reduction at heart, we have $(M_{\lambda})\git G=(M\git
G)_{\lambda}$.  If we express the equivariant volume of the reduced space $M\git G$ as the limit of the equivariant volumes
of the cut spaces $(M\git G)_{\lambda}$, the equivariant residue formula allows us to calculate the equivariant
volumes of $M\git G$ as 

\begin{equation}\label{master}
\begin{split}
\int_{M\git G}e^{\omega +\mu_{T}} & = \lim_{\lambda \rightarrow \infty} \int_{(M\git
G)_{\lambda}}e^{\omega+\mu_{T}}\\
&= \lim_{\lambda \rightarrow \infty} \int_{(M_{\lambda})\git G}e^{\omega+\mu_{T}}\\
&= \lim_{\lambda \rightarrow \infty}\sum_{F\subset \left(M_{\lambda}\right)^{T_{G}}}\frac{K}{\abs{W}}\JKres^{\Lambda}\varpi^{2}
\int_{F}\frac{e^{\omega+\mu_{T}+\mu_{T_{G}}}}{e_{T\times T_{G}}(\nu_{F})}.
\end{split}
\end{equation}
 
 One now has to make the basic but important observation that the symplectic cut introduces new fixed points --- under the conditions
stated above we have, for large $\lambda$, $$M_{\lambda}^{T_{G}}=M^{T_{G}}\coprod \left(M\git_{\lambda} U(1)\right)^{T_G}. $$  Furthermore, when considering the limit as
$\lambda\rightarrow \infty$ we restrict ourselves to the open cone in $\mathfrak{t}_{\C}$ as in section~\ref{pw}.  On this open
cone the limit of most --- but not all --- of the terms given by the residue formula vanishes.

This is somewhat reminiscent of the physical concept of
renormalization for quantum field theories: we have introduced a
`scale' (the momentum cut $\lambda$), but it turns out that there are quantities independent of the scale, and those are the ones that survive the limit as $\lambda$ goes to infinity.

We further remark that for our purposes, we shall use theorem~\ref{equiJK} in the case where the second group $H=T$ is a torus.  In this situation, the integrals over the fixed point components are integrals in $T$-equivariant cohomology that we can calculate via the usual localization.  Therefore we can actually replace the sum on the right-hand side with the sum over all the $T_G\times T$-fixed point components $F\subset M^{T_G\times T}$.

\subsection{Basic Example}\label{basex}
Let us illustrate this principle in the simplest possible case:
the equivariant integral on $N=\C$ with respect to a circle $H=T^1$ acting linearly with weight two.  This is nothing but a Gaussian integral, which can be calculated using the Prato-Wu theorem
$$\int_{N}e^{\omega+\mu_{T}(\tau)}=\int_{\C}e^{-2\nrm{z}^{2}\tau}\frac{i}{2}dz\wedge d\overline{z}=\frac{1}{2\tau},$$where the right-hand side is the contribution of the unique fixed point to the localization formula.
This result needs to be interpret as holding only on the open cone $\tau\in(0,\infty)$ in $Lie(T)\cong \R$, or if preferred on the cone $(0,\infty)\times i\R$ in the complexification $\mathfrak{t}_{\C}$.  Indeed, if we would follow the reasoning of equation~(\ref{cmp}) and compactify $N$ through a symplectic cut to obtain $N_{\lambda}=\CP_1$, then this cut would introduce a single extra fixed point, with contribution $\frac{e^{-\lambda\tau}}{-2\tau}$.   In the limit $\lambda\rightarrow\infty$, the contribution of this fixed point `at infinity' would vanish for $\tau \in (0,\infty)$ or $(0,\infty)\times i\R$.

\begin{figure}[!htb]\begin{center}\leavevmode
\begin{xy}
\WARMprocessMoEPS{volgende}{eps}{bb}
\xyMarkedImport{volgende}
\xyMarkedMathPoints{1-9}
\end{xy}
\caption{\label{voorbeeldje}Moment maps for the groups $G$ and $H$}
\end{center}
\end{figure}
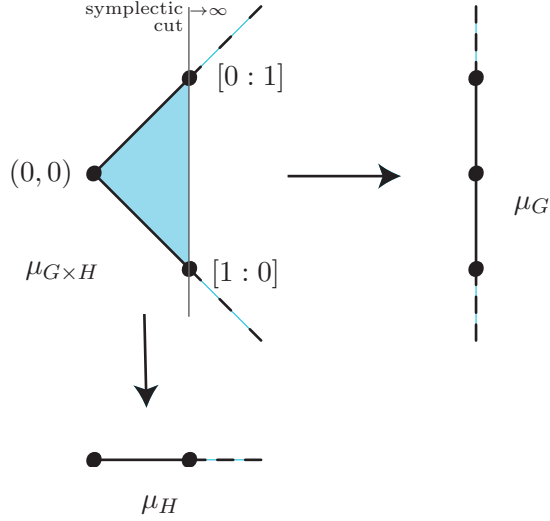

Alternatively, we can construct
$N$ as well as the reduction of $M=\C^{2}$ by the action of a circle $G=T^1$ acting with weights $1$ and $-1$: $$s(z_1,z_2)=(sz_1,s^{-1}z_2).$$  The action of $H$ lifts to $M$ as well, by $$t(z_1,z_2)=(tz_1,tz_2).$$ We can depict the images under the moment maps of $M$ and $M_{\lambda}$ as in figure~\ref{voorbeeldje}.

After making the cut $M_{\lambda}\cong \C^2\coprod \CP_1$ and using affine coordinates for $\C^2$ and projective coordinates for $\CP_1$, it is easy to see that there are three fixed points in $M_{\lambda}$: the original $(0,0)$ in $M$, and the two new fixed points, `at infinity', $[1\!:\!0]$ and $[0\!:\!1]$.  For
the calculation of the residue, we only need to take $[1\!:\!0]$ into account, as the residues of the contributions of the other fixed points give zero.  The weights of the
isotropy representation of $T$ on its tangent space are directly read off from figure~\ref{voorbeeldje}: $\sigma-\tau$ and
$2\sigma$.
 Hence we get as an integral: 
\begin{equation*}
\begin{split}\int_{M\git H} e^{\omega+\mu_{T}} & = \lim_{\lambda\rightarrow \infty}\int_{M_{\lambda}\git H}
e^{\omega+\mu} \\ &=\lim_{\lambda\rightarrow \infty}  \JKres^{\Lambda} \left(\frac{e^{-\lambda(-\sigma+\tau)}}{(\sigma-\tau)(2\sigma)}\right)\\
& = \lim_{\lambda\rightarrow \infty}\left( \frac{e^{-\lambda\tau}}{-2\tau}+\frac{1}{2\tau}\right)\\
 &=\frac{1}{2\tau}.
\end{split}
\end{equation*}
The reason only one of the two terms given by the residue survives the limit as $\lambda\to\infty$ is exactly that the cone on which we consider the outcome (namely $\tau\in(0,\infty)$ or $(0,\infty)\times i\R$) is such that every exponential that still occurs vanishes.

\section{Equivariant volumes of quotients by linear actions}\label{linact}
In order to put the above idea to use, one must have a way of handling the limit of the residues as $\lambda$ goes to infinity.  Therefore we restrict ourselves from now on to symplectic reductions and \HK quotients of vector spaces by linear group actions, where this becomes feasible.  We first need some remarks about how we interpret the symplectic cut for \HK quotients.
\subsection{Complex Moment Maps and \HK quotients}
We will be interested in applying this philosophy to calculate equivariant volumes  of \HK quotients of a vector space $V$ by a linear group action,
$$V\hk_{(\lambda_1,\lambda_2,\lambda_3)} G = \mu_1^{-1}(\lambda_1)\cap \mu_2^{-1}(\lambda_2)\cap \mu_3^{-1}(\lambda_3) / G $$  However, we will not be interested in the full \HK structure, but rather in the symplectic structure induced by one of the complex structures.  With a choice of a preferred complex structure we can single out the moment map $\mu_1$ associated with the corresponding \K form, and pack the `other' moment maps together in a complex valued $\mu_{\C}= \mu_2+i\mu_3: V \rightarrow \mathfrak{g}^*\otimes \C$.  We will therefore mainly think of the \HK quotients as symplectic reductions of a level set of this complex moment map:
$$V\hk_{(\lambda_1, \lambda_2,\lambda_3)} G = \mu_{\C}^{-1}(\lambda_2+i\lambda_3) \git_{\lambda_1} G$$
The torus actions of $T$ on $V$ that we will consider do not preserve the \HK structure, but as $\mu_{\C}$ is (quadratic) homogeneous it will preserve the variety $\mu_{\C}^{-1}(0)\subset V$ and therefore it will descend to an action on $V\hk_{(\lambda_1,0,0)} G= \mu^{-1}_{\C}(0)\git_{\lambda_1} G$.

Now in order to apply the calculational scheme outlined above, we need to study the symplectic cut $\left(\mu^{-1}_{\C}(0)\right)_{\lambda}$.  Our calculations become much simpler, however, if we remark here that this symplectic cut, a subspace of $V_{\lambda}$, can still be interpreted as a zero-set of a section of a vector bundle.  Indeed, in the cases studied later on, the symplectic cut happens with respect to a circle acting with global weight, and hence the cut space $V_{\lambda}$ is a projective space.  As $\mu_{\C}$ is quadratic, we can therefore associate with it an equivariant section of the $T$-equivariant bundle $\mathcal{O}(2)\otimes \mathfrak{g}_{\C}$ over the projective space $V_{\lambda}$, which we will also denote as $\mu_{\C}$: 
\begin{equation}\label{lim} \xymatrix{ \mathfrak{t}^*_{\C}\otimes \mathcal{O}(2) \ar[dd]\\ \ \\ V_{\lambda}\ar@/^1.5pc/[uu] ^{\mu_{\C}} } 
\end{equation}
By further abuse of notation we will say that $$(V_{\lambda})\hk G = \mu^{-1}(0)\git G.$$
This allows us now to reduce integration over $\left(V\hk
  G\right)_{\lambda}$ to integration over $\left(V\git
  G\right)_{\lambda}$, by using the (equivariant) Euler class
of the bundle as equivariant Poincar\'e dual to this zero-set.
Strictly speaking there is no Poincar\'e duality in equivariant cohomology, but one still has the property that $\int_{M} e_G(E) \wedge \alpha = \int_{s^{-1}(0)} i_s^* \alpha$ for a sufficiently general\footnote{For the spaces that we are interested in this is indeed the case: it was shown by Crawley-Boevey~\cite[Theorem 1.1 and remarks to $\S 1$]{craw-boev} that the zero-level of the complex moment map for ADHM spaces is a complete intersection.} section $s$ of an equivariant bundle $E$ over a compact $M$, see e.g.~\cite{meinrequi}.
Hence we can adapt the calculational method given above in~(\ref{master}) to the case of \HK quotients of vector spaces:
\begin{equation}\label{ideetje}
\begin{split}
\int_{V\hk G} e^{\omega+\mu_T} & = \lim_{\lambda\rightarrow \infty} \left( \int_{\left( V \hk G \right)_{\lambda} } e^{\omega+\mu_T}\right) \\
& = \lim_{\lambda\rightarrow \infty} \left( \int_{\left( V_{\lambda} \right) \hk G } e^{\omega+\mu_T} \right) \\
& = \lim_{\lambda\rightarrow \infty} \left( \int_{\left( V_{\lambda} \right) \git G } e^{\omega+\mu_T}\ e_T(\mathfrak{g}_{\C}\otimes \mathcal{O}(2)) \right) 
\\
& =  \lim_{\lambda\rightarrow \infty} \left( \frac{K}{\abs{W} } \sum_{F\subset \left(V_{\lambda}\right)^{T_G\times T}} \JKres^{\Lambda} \ {\scriptstyle{ \varpi^2}} \int_F {\scriptstyle{ \frac{e^{\omega+\mu_T+\mu_{T_G}}\ e_{T_G\times T}(\mathfrak{g}_{\C}\otimes \mathcal{O}(2))}{e_{T_G\times T} (\nu_F)}}} \right)\!\!.
\end{split}
\end{equation}

We have been conspicuously silent about an important aspect of
this method: the possibility that the quotient $M\git G$ or
$V\hk G$ is singular. (The simple example discussed in
section~\ref{basex} above was a particular case where a
reduction at a singular value of the moment map was
nevertheless smooth.)  In essence, everything still goes through, and we postpone the discussion about this to section~\ref{nstntn}.

\subsection{Calculational strategies}\label{strat}

Since in practice most, if not all, finite-dimensional \HK quotients are quotients of vector spaces by linear group actions, we will now restrict ourself to those cases. For these, as outlined above, we need a circle inside $T$ that acts with a single global weight, and hence the cut space $V_{\lambda}$ will be a projective space $V_{\lambda}=V\coprod \proj V$.  As all the fixed point data  in this projective space are determined by the representation of $T_G\times T$ on $V$, we  can expect a more direct approach that bypasses the procedure of making the cut to compactify,  examining  the fixed points on the cut locus $\proj V$ and their residues, and then taking the limit $\lambda\rightarrow \infty$.  Indeed this is the case, and below we give a description for the outcome 
that no longer mentions the cut. We will take a shamelessly
pragmatic approach, explicitly using coordinates and working
with a residue operation defined through iterated one-dimensional residues as in (\ref{beupf}).  For simplicity we will assume that all relevant (i.e. whose residue is not manifestly zero) $T\times T_G$ fixed points on $V_{\lambda}$ are isolated --- which is the case for all the examples we will discuss below.

We will use the following notation: for a meromorphic function $\frac{p}{\prod \rho_i}$ expressed in linear coordinates $(X_1,\ldots,X_n)$ and with all $\rho_i$ linear, denote by $$\mres^{X_i}_j\frac{p}{\prod \rho_i}$$  the meromorphic function in the variables $(X_1,\ldots,X_{i-1},X_{i+1},\ldots,X_n)$ obtained by taking the residue in the variable $X_i$ of the pole determined by the $j$-th factor in the denominator (and meanwhile considering all other variables as generic constants).  The new form lacks a $j$-th factor in the denominator, but for convenience we keep the same name and labelling of factors (i.e. the index skips $j$).  As an example we have  \begin{multline*}\mres^{X_2}_{3} \frac{X_1-4X_2}{(-X_1+X_2+X_3)(X_1-X_2+X_3)(X_1+X_2-X_3)}\ \\ =\ \frac{X_1 - 4(X_3-X_1)}{(-X_1+(X_3-X_1)+X_3)(X_1-(X_3-X_1)+X_3)}.\end{multline*}

\subsubsection{Fixed points in $V_{\lambda}$} Now, let us
examine how to implement the above calculational method
(\ref{ideetje}) when we have a (compact connected) Lie group
$G$ with maximal torus $T_G$ which acts tri-Hamiltonianly on a
\HK vector space $V$, with complex moment map
$\mu_{\C}$. (Everything below still goes through if $V$ is just
a complex vector space of which we want to take a symplectic
reduction, by just dropping the extra factors coming from
$\mu_{\C}$ throughout.)  We assume that $V$ is further equipped
with the action of a torus $T$, commuting with $G$, such that
$T$ contains a circle acting with global weight one.  Let
$\rho_i$ be the weights (repeated in the case of
multiplicities) occurring in the weight decomposition of the torus action of $T\times T_G$ on $V$.  Clearly the `new' fixed point components of the cut space in $V\git_{\lambda} U(1) \subset V_{\lambda}$ will correspond exactly to weight spaces $V_{\rho_i}$ for some $i$.
With the assumption of isolated fixed points made above, the relevant such weights will have multiplicity one, and the image under the moment map $\mu_{T\times T_G}$ of such an isolated fixed point of $V_{\lambda}$ will be $\lambda \rho_i$.  One can also easily see that the weights of the isotropy representation at the fixed point corresponding to the weight $\rho_i$ are given by $\rho^i_1=\rho_1-\rho_i,\dots, \rho^i_{i-1}=\rho_{i-1} - \rho_i, \rho^i_i=-\rho_i,\rho^i_{i+1}=\rho_{i+1}-\rho_i,\dots$.  We will further abuse notation by writing the equivariant Euler class of the bundle (\ref{lim}), of which the complex moment map is a section, at the fixed point $0$ as $e(\mu_{\C})$ as well (this is also a polynomial, a product of linear factors), and at the fixed point determined by $\rho_i$ as $e(\mu_{\C}^i)$. 

\subsubsection{Admissible paths} 
\begin{figure}[!htb]\begin{center}\leavevmode
\begin{xy}
\WARMprocessMoEPS{admpath2}{eps}{bb}
\xyMarkedImport{admpath2}
\xyMarkedMathPoints{1-9}
\end{xy}
\caption{\label{admispath}Example of a path $P$, taking successive residues at various poles  }\end{center}
\end{figure}
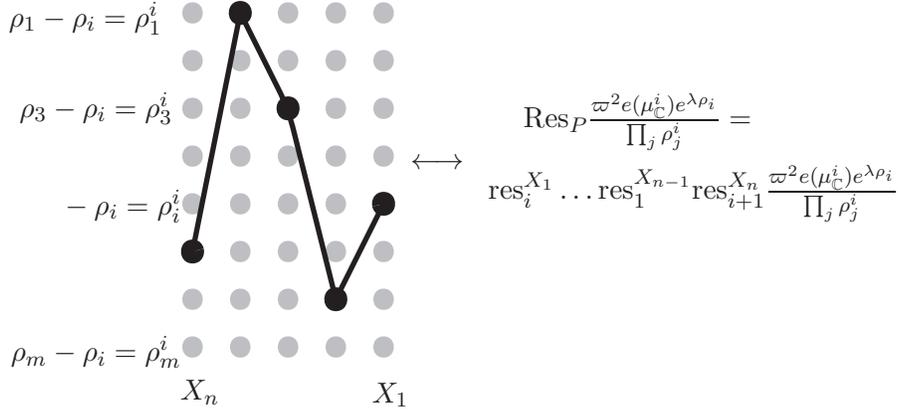

Let us examine the contributions of such a fixed-point component to the formula given above.  Clearly the iterated residue $\JKres^{\Lambda}$ (\ref{beupf}) gives rise to a sum of terms, each of which is obtained by taking the residue at one pole for each variable $X_i$ on $\mathfrak{t}^*$:  
\begin{equation}\label{sommeke}\JKres^{\Lambda} \frac{\varpi^2e(\mu_{\C}^i) e^{\lambda\rho_i}}{\prod_j \rho^i_j} =  \sum \mres^{X_1}_{j_1}\dots \mres^{X_n}_{j_n} \frac{\varpi^2e(\mu_{\C}^i) e^{\lambda \rho_i} }{\prod \rho^i_j}.\end{equation}  Calculating (\ref{sommeke}) in practice comes down to determining exactly which iterated residues need to be included in the right-hand side.  For each term in the right-hand side of (\ref{sommeke}) we can indicate which poles were chosen for which variables by a path $P$ as in figure~\ref{admispath}, and denote this by $\res_P$.

As taking the residue at a pole eliminates that pole, each path can only contain every pole at most once, but certainly this is not a sufficient condition.  
Let us examine more closely which paths can occur in the sum (\ref{sommeke}).
Recall from (\ref{bijna2}) the definition of the
one-dimensional residue $$\jkres^+(f(x) dx)=\sum_{\lambda_j
  \geq 0} \sum_{b\in \C} \mres_b (g_j(x) e^{\lambda_j x}).$$
Iterating this gives the condition that just before taking the
residue for the $X_j$ variable determined by the path, the
coefficient of the $X_{j}$ variable of the exponent has to be
positive, and this holds for all $j$.  Let us call this condition A.

Furthermore, we are only interested in the terms that survive
the limit $\lambda\rightarrow \infty$ (as always, only
considering the answer as a function on an open cone in
$\mathfrak{t}^*$).  In order for there to be any term that
survives the limit, the last pole that is evaluated has to be
the $i$-th one, so our path has to end at the $i$-th
vertex in the right-most column. Call this condition B. Indeed, if the $i$-th pole
occurs sooner in the path, the next step in the full iterated
residue $\JKres^{\Lambda}$ clearly yields zero.  Likewise, if
the $i$-th pole is never used in a residue coming from  a path,
the term corresponding to that path would still have a non-zero
coefficient in the exponent, and the cone in the Lie algebra on which we evaluate the final result is such that in the limit $\lambda\rightarrow \infty$ the term would vanish. 

Let us call a path satisfying both these conditions A and B an \emph{admissible path}.  Hence we have
 \begin{equation*}\lim_{\lambda \rightarrow \infty} \JKres^{\Lambda} \frac{\varpi^2 e(\mu_{\C}^i) e^{\lambda\rho_i} }{\prod_j \rho^i_j}\ =\  \sum_{\text{\parbox{1in}{\begin{center}admissible paths~$P$ ending at $\rho_i$\end{center}}}} \res_P \frac{\varpi^2e(\mu^i_{\C}) e^{\lambda\rho_i}}{\prod_j \rho^i_j}\end{equation*}

\subsubsection{Yoga of admissible paths}  From the above it is clear that in practice calculating the residue $\JKres^{\Lambda}$ boils down to determining which paths $P$ are admissible, and calculating the corresponding contribution $ \res_P \frac{\varpi^2e(\mu^i_{\C})e^{\lambda\rho_i}}{\prod_j \rho^i_j}$ for each admissible path $P$ ending at $\rho_i^i$.   Let us use the following notation: given a linear form $\rho$ and a choice of (linear) coordinates $(X_1,\ldots,X_n)$ on a vector space, we say that $\rho>_j 0$ if the $X_j$ coefficient of $\rho$ is greater than $0$.
With this notation, we can state the conditions for a path to
be admissible as follows: first, we need to verify that
$\rho_i>_n0$.  After taking the residue in the $X_n$ variable
at the pole determined by $\rho^i_{j_1}$, we then want $\rho_i>_{n-1}0$, and this keeps repeating until variable $X_{n-1}$.  In particular, when testing if a path is admissible, this requires at each step calculating $\mres^{X_l}_{j_l}$ before we can determine if the vertex $(l, j_l)$ can be part of the path, which is quite inefficient.  

Using some elementary linear algebra, however, we can reshuffle the calculations for contributions of admissible paths in such a way that we gain two distinct advantages.  Firstly, the conditions for a path to be admissible become more straightforward.  Secondly, we can apply the iterated residues corresponding to all admissible paths \emph{to the same single function}, irrespective of the vertex at which the paths end (i.e.\ the fixed point to whose contribution the term belongs).  Hence we eliminate referring to the actual `fixed points at infinity', as we directly extract their contributions from one single function.

\begin{lemma}
A path $P=\{ (X_n,j_n), \ldots (X_2,j_{2}), (X_1,i)\}$ giving rise to the term
$$\res_P\frac{\varpi^2e(\mu^i_{\C})e^{\lambda\rho_i}}{\prod_j \rho^i_i}=\mres^{X_1}_ {i}\mres^{X_{2}}_{j_{2}}\dots \mres^{X_n}_{j_n}\frac{\varpi^2e(\mu^i_{\C})e^{\lambda\rho_i}}{\prod_j \rho^i_j}$$
in the contribution of the fixed point corresponding to the weight $\rho_i$ is admissible if and only if we have at the beginning $\rho_i>_n 0$ and then at each step $\rho_{j_l}>_{l-1} 0$.  Furthermore, its contribution can also be computed as
\begin{equation*}\begin{split}\res_P\frac{\varpi^2e(\mu^i_{\C})e^{\lambda\rho_i}}{\prod_j \rho^i_i}&=(-1)^n \mres^{X_1}_{j_{2}}\dots\mres^{X_{n-1}}_{j_n}\mres^{X_n}_i\frac{\varpi^2 e(\mu_{\C})}{\prod_j \rho_j}\\ &=(-1)^n \res_{\widetilde{P}}\frac{\varpi^2e(\mu_{\C}) }{\prod_j\rho_j}\end{split}\end{equation*}
where $\widetilde{P}$ is the path $\{ (X_n, i), (X_{n-1}, j_n), \ldots , (X_1, j_2)\}$.
\end{lemma}
\begin{proof}  For the first statement it suffices to remark
  the following.  We each time want to test the function $\rho_i$ in the exponent after consecutive residues: $\rho_i>_n 0, \rho_i>_{n-1} 0, \dots$.  However the poles at which we take the consecutive residues are exactly $\rho^i_{j_n}=\rho_{j_n}-\rho_i,\rho^i_{j_{n-1}}=\rho_{j_{n-1}}-\rho_i,\dots$.  As taking a residue at such a pole essentially sets that term to zero, we might as well test for $\rho_i>_n 0, \rho_{j_n}>_{n-1}0, \rho_{j_{n-1}}>_{n-2} 0,\dots$.

For the second statement, notice that calculating one-dimensional residues $\mres^{X_l}_{\rho^i_j}$ can essentially be understood in terms of elementary linear algebra as Gaussian elimination for a matrix.  Indeed, if you organize all the occurring linear functions (expressed in coordinates) as rows in a matrix as follows:
$$
\begin{array}{cc}
	\   & \!\!\!\!\!\!\!\!
	\begin{array}{ccc}
		\overbrace{
			\begin{array}{ccc} 
				X_n & \ldots & X_1 
			\end{array}
		}^{\text{\parbox{1.3in}{\tiny variables used for residues}}} &  \!\!\! \!\!\!
		\overbrace{
			\begin{array}{ccc}
				Y_1 & \dots &  Y_m 
			\end{array}
		}^{\text{\parbox{1.3in}{\tiny \ \ \ \ \ remaining variables}}}
		 & \
	\end{array} \\
	\begin{array}{c}  \!\!\!\!\!\!\!\text{\parbox{.4in}{\tiny factors used for poles} } \
		\begin{cases} 
			\rho^i_{j_n}\ & \  
			\\ \vdots \ & \ 
			\\ \rho^i_{j_{2}} \ & \
			\\ \rho^i_i \ & \
		\end{cases}\\ 
		 \!\!\!\!\!\!\!\text{\parbox{.4in}{\tiny all other factors} } \
		\begin{cases} \phantom{\rho^i_{j_{2}}}
			\ & \  
			\\  \ & \ 
		\end{cases}
	\end{array} \!\!\!\!\!\!\!\!\!\!\!\!\!\!\!\!\!\!\!\!\!\!\! & 
	\begin{array}{cc}
		\left[\begin{array}{cccccc}
			\phantom{X_1} \ \ \ \ & \phantom{\dots}  & \phantom{X_n}  & \phantom{Y_1}  &\phantom{\dots}  &\phantom{Y_m}  \\ 
			\ & \ & \ & \ & \  &\ \\   
			\ & \ & \ & \ & \  &\ \\   
			\ & \ & \ & \ & \  &\ \\   
			\ & \ & \ & \  & \ &\ \\  
			\ & \ & \ & \ & \ & \ \\
			\ & \ & \ & \ & \ & \ 
		\end{array}\right] & \!\!\!\!\!\!
		\begin{array}{l}
			\ \phantom{\rho^i_{j_1} }
			\\ 
			\ \phantom{\vdots }
			\\
			\ \phantom{\rho^i_{j_{n-1}} }
			\\ \ \phantom{\rho^i_i }
			\\ \ \\ \
 		\end{array}
	\end{array}
\end{array}
$$
 then taking a residue in a variable can be interpreted as doing the column-e\-li\-mi\-na\-tion for the corresponding column (given by the variable used for the residue) and row (given by the factor), and replacing all but the used row in the original function with the corresponding new rows.  The factor of the pivotal entry is all that remains of the row used.
   
From this one immediately sees that taking the consecutive one-di\-men\-sio\-nal residues comes down to doing row-reduction for the first $n$ rows, replacing the remaining factors in the original function by the corresponding new rows, and dividing the function by the determinant of the upper $n \times n$ submatrix in the above matrix.  
With this in mind it is now clear that one can in fact change
the ordering of the poles used for the consecutive residues, as
well as the ordering of the variables used, with  the sign
change that occurs through the determinant as the only
penalty.  The single condition that needs to be satisfied is that for the new ordering, all the pivotal elements need to be non-zero.

Now, given our ordering of the poles (i.e. $\rho^i_{j_n}, \dots, \rho^i_{j_{2}}, \rho^i_i$), change the ordering by making the last pole used first (hence obtaining $\rho^i_i,\rho_{j_n}^i,\dots, \rho^i_{j_2}$).   As we have used $\rho^i_i=-\rho_i$ for the first pole, this immediately implies that we might as well change all the other factors (i.e. $\rho^i_{j_n}, \dots, \rho^i_{j_{2}}$ as well as $e(\mu^i_{\C})$) to the corresponding factors (i.e. $\rho_{j_n},\dots,\rho_{j_2}$ and $e(\mu_{\C} )$) for the `central function' $\frac{\varpi^2 e(\mu_{\C})}{\prod_j \rho_j}$.
This ordering also has the advantage that it exactly corresponds to the ordering of the testing explained above, and furthermore the consecutive conditions that the pivotal elements are non-zero is indeed implied through the consecutive tests $\rho_{j_l}>_{l+1}0$.
\end{proof}
With this it makes sense to define the following:
\begin{definition}
Let $X_1,\ldots,X_n$ be a choice of linear coordinates on a
complex vector space $V$.  Then for any rational function $f=\frac{p}{q}$ on $V$,  where $q=\prod_j \rho_j$, and where each $\rho_j$ comes with a preferred polarization, we define $\res_+^{X_1,\ldots X_n}$ inductively as $$\res_+^{X_1,\ldots X_n}= \mres_+^{X_1}\dots \mres_+^{X_n} \frac{p}{q}$$ where in each step the other variables are treated as generic constants, and where the single-variable residues are defined as $$\mres_+^X \frac{p}{q}=\sum_{\rho_j >_X 0} \mres^X_{j} \frac{p}{q},$$ $\rho_j$ in the sum on the right hand side being the factors of $q$.%, and  $\mres^X_{\rho_j}$ being the residue in the variable $X$ for the pole of $\frac{p}{q}$ determined by the factor $\rho_j$.
\end{definition}
With this definition we now can summarize the discussion of
this section: 
\begin{theorem}Let a compact connected Lie group $G$ act tri-Hamiltonianly on a \HK vector space $V$, such that the group action commutes with the action of a torus $T$.  Assume that $T$ acts Hamiltonianly with respect to one of the \K structures on $V$.  Let $\omega$ and $\mu_T$ denote the corresponding \K form and moment map  on the \HK quotient $V\hk G$.  Then with the notation as above we have $$\int_{V\hk_{(0,0,0)} G} e^{\omega+\mu_T}=\frac{\tilde{K}}{\abs{W_G}}\res_+^{X_1,\ldots, X_n} \frac{\varpi^2 e(\mu_{\C})}{\prod_j \rho_j},$$ 
where $\tilde{K}$ is the constant $\frac{(-1)^{n_++n}}{\volu(T_G)}$.\end{theorem}
Here $X_1,\dots,X_n$ is a linear coordinate system on the Lie algebra of $G$, suitable in the sense of (\ref{beupf}).  For the volume of  a symplectic reduction of a vector space with a linear action, the formula is identical except that $e(\mu_{\C})$ is missing.

We remark here the difference in the definitions of the various residue operations:  $\JKres^{\Lambda}$ and the one-dimensional $\jkres^+_X$ on the one side, and $\res_+^{X_1,\ldots X_n}$ on the other. The former, originally coming from a Fourier-transform, involve a positivity test for the coefficient in the exponent of the numerator (hence the choice of the position of $\Lambda$ and $+$ in the notation). The latter involves a test of positivity on the factors of the denominator.

The advantage of the reordering of the poles is that in practice it offers an easier way to go through the calculations of residues: at each step (i.e.\ when taking the residue for each successive variable) just take the sum of the residues at the poles corresponding to factors where the variable has a positive coefficient.
We remark that the way we have written the residue above does not simply apply to rational functions --- in particular each factor in the denominator has to come with a preferred sign, which is similar to the choice of polarization used in the Guillemin-Lerman-Sternberg approach to the Duistermaat-Heckman formula~\cite{symfib}.  In the next section we will use this result further to compare our calculational method with the one given by Nekrasov and Shadchin~\cite{ABCD} in the physics literature.

\section{Instanton Counting}\label{nstntn} We are now ready to apply the calculational techniques described above to 
our guiding examples: the calculation of the equivariant volumes of the ADHM spaces as they occur in instanton counting.  We discuss the corresponding residue formulas, both for their own sake and as examples of how to implement the calculational heuristic given above in various settings, and rederive the results of Nekrasov and Shadchin~\cite{ABCD}.  The most important aspect still requiring settlement is the matter of the singularities of the quotients.  

\subsection{Physical background}
The recent work regarding instanton counting~\cite{nekra} is a direct answer to an open problem in physics.  In particular, it is concerned with $N=2$ supersymmetric quantum Yang-Mills theory. 
One now wants to describe the low-energy behavior of this theory, and using the ideas of the Wilson renormalization group, one finds that this low-energy limit is described by means of another action, the \emph{effective action}.  In general it is extremely difficult to explicitly write down the effective Lagrangian.  
However, in the case of $N=2$ quantum Yang-Mills theory it was shown by Seiberg and Witten in the seminal~\cite{SW} that the entire Lagrangian for the low energy effective theory is determined by a single (multivalued) holomorphic function $\mathcal{F}_{SW}$, the Seiberg-Witten prepotential.  

Furthermore, using an assumption that a form of electric-magnetic duality (called Montonen-Olive duality) holds for the quantum theory, Seiberg and Witten managed to describe this prepotential analytically, in terms of period integrals of a family of curves.  Their work was a major breakthrough in physics, but led to spectacular advances in mathematics as well, in particular in the framework of  the Seiberg-Witten invariants in differential geometry.

The original work of Seiberg and Witten was based on an assumption of duality, and though this duality is widely believed to hold, it is conjectural even for physicists as no derivation from first principles is known. In the last 10 years attempts were made to derive and verify their results directly.  It was shown that this reduces to calculating certain integrals over moduli spaces of framed instantons on $\R^4$ (see e.g.~\cite{calc-inst} for background), but the actual calculations of these integrals were very difficult in general.

The solution to this outstanding problem was finally accomplished by Nekrasov in 2002~\cite{nekra}.  The strategy Nekrasov employs is to use maximal symmetry on the moduli spaces, induced by change of framing of instantons and rotations in $\R^4$, and then compute equivariant volumes with respect to this group action by means of localization techniques in equivariant cohomology.  

In particular Nekrasov considers the following generating function\footnote{The full generating function Nekrasov studies  $Z=Z^{\text{pert}}Z^{\text{inst}}$ also has an extra factor, $Z^{\text{pert}}$, the \emph{perturbative one-loop contribution}, which is of no concern to us.}
\begin{equation}
\label{nekrapart}
Z^{\text{inst}}(q,\tau,\epsilon_1,\epsilon_2)=\sum_{k=0}^{\infty} q^k \int_{\modu{n}{c}} e^{\omega+\langle\mu_{T},(\epsilon_1,\epsilon_2,\tau)\rangle},\end{equation}
where $\modu{n}{c}$ is the instanton-moduli space of rank $n$,
charge $c$ framed instantons on $\R^4$.  The torus $T=T^2\times
T_G$ is the product of $T^2$,  the maximal torus of $SU(2)$
that acts diagonally on $\R^4$ after the identification
$\R^4=\C^2$ and hence has an induced action on $\modu{n}{c}$,
and $T_G$, the maximal torus of the gauge group that acts on
$\modu{n}{c}$ by changing the framing of the instantons at infinity.  Here $\epsilon_i$ and $\tau$ are coordinates on the Lie algebras of $T^2$ and $T_G$ respectively.  With the interpretation given in section \ref{pw} above this is a mathematically well-defined object, a function on an open subset of the Lie algebra of $T_2\times T_G$, which by analytic continuation gives a meromorphic (even rational) function on the whole of the complexified Lie algebra of $T_2\times T_{G}$.

Nekrasov argues that one can write $Z^{\text{inst}}(q,\tau,\epsilon_1,\epsilon_2)$ as $$Z^{\text{inst}}(q,\tau,\epsilon_1,\epsilon_2)=\exp\left(\frac{\mathcal{F}^{\text{inst}}(q,\tau,\epsilon_1,\epsilon_2)}{\epsilon_1\epsilon_2}\right)$$ where the function $\mathcal{F}^{\text{inst}}(q,\tau,\epsilon_1,\epsilon_2)$ is analytic and regular near $\epsilon_1,\epsilon_2=0$.  Furthermore, he claims that $\mathcal{F}^{\text{inst}}|_{\epsilon_1,\epsilon_2=0}$ corresponds to the instanton part of the prepotential of Seiberg-Witten.  As the latter can be defined rigorously in terms of periods of certain families of curves, this correspondence gives rise to a remarkable conjecture in geometry.  In 2003, this conjecture was proven independently by Nakajima and Yoshioka ~\cite{naka,naka2} and Nekrasov and Okounkov~\cite{rndprt} for the gauge group $SU(n)$, using very different techniques.  In~\cite{ABCD} Nekrasov and Shadchin indicated how the proof of~\cite{rndprt} could be adapted to the other classical gauge groups.  More recently, a non-computational proof for all simple gauge groups was given by Braverman and Etingof in ~\cite{braver1,braver2}.  

\subsection{Equivariant volumes of ADHM spaces for $SU(n)$ instantons}
\subsubsection{ADHM spaces}  The space \modu{n}{c} we are considering here is the moduli-space of \emph{framed $SU(n)$ instantons} with instanton number $c$ (cf.~\cite{ADHM,naka2}) , constructed as a \HK quotient of
the vector space \adhm{n}{c} of \emph{ADHM data}  $$ \adhm{n}{c} = \{(A,B,i,j)\}$$ 
where $A,B,i,j$ are linear maps between complex hermitian vector spaces $V$ and $W$ of dimension $c$ and $n$ respectively  as follows:
$$\xymatrix{V \ar@(ul,dl)[]_{A}\ar@(ur,dr)[]^{B}\ar@/_/[dd]_{j} \\  \\ W\ar@/_/[uu]_{i} \\}$$
and the quotient is taken by the group action of $U(V)$ which acts on $A$ and $B$ by conjugation and on $i$ and $j$ by left and right multiplication.

Strictly speaking, the moduli space of framed instantons (or,
equivalently, framed rank $n$ holomorphic bundles on $\CP_2$
with second Chern number $c$) is the non-singular locus of this
quotient.  This non-singular locus doesn't satisfy the conditions
to apply either the Prato-Wu theorem (\ref{pwform}) or the
method described above --- in particular the moment map for the
torus actions we are considering is not proper, due to the fact that the space is not metrically complete in the \K  metric.    

One can however extend (`partially compactify') these spaces in various ways: first of all, one can allow for so-called \emph{ideal instantons}, which are interpreted in differential geometry as being (framed) connections whose curvature is concentrated at certain points.  This leads to the \emph{Uhlenbeck space}, which is the full \HK quotient: $$\uhl{n}{c}=\adhm{n}{c}\hk_{(0,0)} U(V) =\mu_{C}^{-1}(0)\git_0 U(V)= \mu_{\R}^{-1}(0) \cap \mu_{\C}^{-1}(0) / U(V).$$ 
 The Uhlenbeck space is highly singular, however.  A better option is to extend the non-singular locus in another way, which provides a desingularization of the Uhlenbeck space.  As the group $U(V)$ occurring in the \HK quotient has characters, we can vary the symplectic reduction or GIT quotient (cf. \cite{gitflips,dolgagitflips}) by changing the value of the real moment map.  This gives the \emph{Gieseker space}
 $$\gies{n}{c} = \adhm{n}{c}\hk_{(+,0)} U(V) =\mu_{\C}^{-1}(0)\git_+ U(V)= \mu_{\R}^{-1}(+) \cap \mu_{\C}^{-1}(0) \git U(V).$$
 The Gieseker space has a modular interpretation in algebraic
 geometry as the moduli-space of framed torsion-free sheaves on
 $\CP_2$, i.e. torsion-free sheaves that come with a fixed
 trivialization on a line $\ell_{\infty}\subset\CP_2$ (in
 physics this space is also thought of as the moduli space of instantons on a non-commutative $\R^4$ \cite{noncominst}).  Of particular relevance for us is that $\gies{n}{c}$ is smooth. 
 
\subsubsection{Torus actions on ADHM spaces} 
In instanton counting one is interested in the group action given by $T^2\times T_{U(n)}$, where $T_{U(n)}$ changes
the framing at infinity by acting by the maximal torus of the gauge group (for notational convenience it is easier to
allow this to be $U(n)$ rather than $SU(n)$), and $T^2$ acts by rescaling $\C^{2}\subset \CP^{2}$.  In order to lift this action to the space of ADHM data, we must examine the monad constructed out of ADHM data (see e.g.~\cite[Chapter 2]{nakabook}\cite{instgit}).  The monad is a sequence of vector bundle maps over $\CP_2$:
\begin{equation*}
V\otimes \mathcal{O}(-1) 
\xrightarrow{ {{\displaystyle a=}{
\left(\!\! {\renewcommand{\arraystretch}{0.4}
\begin{array}{c} \scriptscriptstyle x_0 A-x_1 \\ \scriptscriptstyle x_0 B-x_2 \\ \scriptscriptstyle x_0 j 
\end{array} }\! \! 
\right)
}}}
{\renewcommand{\arraystretch}{0.7}
\begin{array}{c}
V\otimes \mathcal{O} \\ \oplus  \\ V\otimes \mathcal{O}\\ \oplus \\ W\otimes \mathcal{O} 
\end{array}}
\xrightarrow{{\displaystyle b=}
\left(\! \! {\renewcommand{\arraystretch}{0.4}
\begin{array}{ccc} \scriptscriptstyle -x_0 B+x_2 & \scriptscriptstyle x_0 A-x_1 & \scriptscriptstyle x_0 i 
\end{array} } \! 
\right)
 }
 V\otimes \mathcal{O}(1).
\end{equation*}
If a set of ADHM data satisfies the vanishing of the complex moment map, this sequence of bundles on $\CP_2$ is actually a complex, and the cohomology  $E=\ker b/\im a$ is a bundle (or sheaf) on $\CP_2$, which indeed satisfies the required triviality on the line at infinity.

By examining the
effect of $T^{2}$ on the monad constructed out of the ADHM data, and by interpreting $W\cong
H^{0}(\ell_{\infty},E|_{{\ell_{\infty}}})$
as the trivialization on the line at infinity $\ell_{\infty}$ 
by means of the Beilinson spectral
sequence
(see e.g.\cite{nakabook,okonek}), one can lift the action of $T_{{U(n)}}\times T^{2}$ to \adhm{n}{c} as was done in \cite{naka}:
for $(e_{1},e_{2})\in T^{2}, t\in T_{U_{n}}$ this gives
\begin{equation}\label{goodjobnaka}
(e_{1},e_{2},t).(A,B,i,j) = (e_{1}A,e_{2}B, i t^{-1}, e_{1} e_{2} t j).
\end{equation}

\subsubsection{Volumes}
For this torus action one can now calculate the regularized or equivariant volume of the moduli space.  We remark that the original question asks for the equivariant volume of the moduli space with respect to the \K form it inherits from being included in the Uhlenbeck space (which is an affine variety).  Nevertheless, we can work with the Gieseker space by pulling back the symplectic form from the Uhlenbeck space to the Gieseker space --- so we desingularize in algebraic geometry but not in symplectic geometry.  This gives a closed 2-form on $\gies{n}{c}$ which is degenerate on the exceptional set of $\gies{n}{c}\rightarrow \uhl{n}{c}$;  on $\modu{n}{c}\subset\gies{n}{c}$ it is exactly the form we are concerned with.  Alternatively, one could think of the degenerate symplectic form as a limit of proper symplectic (even \K\!\!) forms that degenerate in the limit, and the regularized volume that we are interested in as the limit of the corresponding equivariant volumes for the Gieseker space.  Despite the degeneracy, one can speak of moment maps with respect to this form.  As usual, the sum of the $2$-form and the moment map determines a cohomology class in the Cartan model of equivariant cohomology, and one can look at the localization formula for the formal exponential of this class, as was already remarked, even with the degeneracy, in~\cite{atibott}.  From our viewpoint, we look at the integral of a function with respect to a volume form, and as $\gies{n}{c}\setminus\modu{c}{n}$ has measure zero we have $$\int_{\modu{n}{c}}e^{\omega+\mu} = \int_{\gies{n}{c}} e^{\omega+\mu}.$$ 
 
For the torus action described above, one can now give explicit
descriptions of the fixed points (which are all isolated), and
the isotropy re\-pre\-sen\-ta\-tions on their tangent spaces,
as was done for $\gies{n}{c}$ in e.g.~\cite{naka}.  In
particular, they can be described very nicely by $n$-tuples of
partitions $\vec{Y}=(Y_1,\ldots,Y_n)$ such that $\sum_i
\abs{Y_i}=k$.    As all these fixed points lie above the same
unique fixed point in $\uhl{n}{c}$ under the desingularization $\gies{n}{c}\rightarrow \uhl{n}{c}$, they all take the same value (i.e.\ $0$) under the moment map for the degenerate symplectic form on $\gies{n}{c}$, and hence we have by the Prato-Wu theorem \begin{equation}\label{gewloc}\int_{\gies{n}{c}} e^{\omega+\mu} = \sum_{\vec{Y}} \frac{e^{\mu(\vec{Y})}}{e_T(\nu_{\vec{Y}})} = \sum_{\vec{Y}} \frac{1}{e_T(\nu_{\vec{Y}})} = \int_{\gies{n}{c}} 1,\end{equation}
where the right-hand side has to be interpreted as a formal integral (of the class $1$) defined by a localization formula.  This is the viewpoint taken by several authors: the integrals that form the coefficients of the Nekrasov partition function are the formal equivariant integrals of $1$ over the Gieseker space $\gies{n}{c}$.  

On the other hand, we could apply the technique described in the sections above.  In order to do this, however, it is very crucial that one thinks of the integral as an equivariant volume rather than the formal integral of 1, as one now has to use the new fixed points `at infinity' introduced by the cut.  The value of these new fixed points under the moment map will not be zero, and through the residue the geometry of the moment map gives the recipe for obtaining their contribution.

From (\ref{goodjobnaka}) we can see that the necessary condition for our cutting con\-struc\-tion to work is clearly satisfied: there is a subgroup in $T^{2}\times
T_{U(n)}$ acting with global weight $1$ on \adhm{n}{c}.  With this lift we can now implement the calculational method
described above.  The resulting formula is
$$\int_{\modu{n}{c}} e^{\omega+\mu}=\frac{1}{c!}\res_+^{\sigma_i}\frac{\varpi^2\ e_T(\mu_{\C})}{e_T(\adhm{n}{c})},$$
where $$\varpi^2=\prod_{e\neq f}\left(\sigma_e-\sigma_f\right),$$ $$e_T(\mu_{\C})=\prod_{1\leq g, h\leq c} \left( \epsilon_1+\epsilon_2+ \sigma_g-\sigma_h\right),$$ 
and 
$$e_T(\adhm{n}{c})= \!\!\!\!\! \prod_{1\leq i, j\leq c} \!\!\!\!\! \left(\epsilon_1\!+\!\sigma_i\!-\!\sigma_j\right) \!\!\!\!\! \prod_{1\leq k,l\leq c} \!\!\!\!\! \left( \epsilon_2\!+\!\sigma_k\!-\!\sigma_l\right) \!\!\!\!\!\!\!\!\!  \prod_{\parbox{.6in}{\begin{center}$\scriptstyle 1\leq m\leq c$ $\scriptstyle 1\leq o \leq n$\end{center}}} \!\!\!\!\!\!\!\! \left(\sigma_m \!-\!\tau_o\right) \!\!\!\!\!\!\!\! \prod_{  \parbox{.6in}{\begin{center}$\scriptstyle 1\leq p\leq c$ $\scriptstyle 1\leq q \leq n$\end{center}}} \!\!\!\!\!\!\!\!\left(\epsilon_1\!+\!\epsilon_2\! -\!\sigma_p\! +\!\tau_q\right) $$
where we use a maximal torus of the form ${\text{diag}(s_1,\dots,s_c)=(e^{\sigma_1},\dots,e^{\sigma_c})}$ for $U(c)$.
\subsection{Equivariant volumes for the other classical gauge groups}
We can now try to implement the same method to compute the regularized volumes of the moduli spaces of $SO(n)$ and $Sp(n)$ instantons.  We remark that for these other classical groups no equivalent of the Gieseker space is known; hence a direct localization calculation as (\ref{gewloc}) is not available.  The Uhlenbeck spaces do exist here, but as before they are highly singular.  We repeat that the integrals we are interested in are the indefinite integrals of a function, $e^{\mu_T}$, depending on a parameter in $\mathfrak{t}_{\C}$ with respect to a volume form $\frac{\omega^n}{n!}$ on the non-singular locus of these Uhlenbeck spaces.  When implementing our method, two minor issues need to be addressed: we need to treat the singularities differently because of this lack of a Gieseker space, and the lifting of the torus action to the space of ADHM data is not automatic.

The ADHM construction for gauge groups $Sp(n)$ and $SO(n)$ was discussed in~\cite{instgit} and described in greater
detail in e.g.~\cite{brysan}.
In both of these cases the groups occurring by which one has to quotient are simple, and hence one cannot hope to obtain a
desingularization of the \HK quotient by means of a variation of GIT quotient as was the case for the gauge group
$SU(n)$.  However, in~\cite{francesdesing}, Kirwan describes a method for constructing desingularizations of singular
quotients $M\git
G$ by blowing up certain subvarieties in $M$ to obtain a new space $\widetilde{M}$, and then constructing the desingularization of $M\git G$ as 
$\widetilde{M}\git G$.  One could therefore use this approach to 
obtain an equivariant desingularization of the Uhlenbeck spaces for gauge groups $Sp(n)$ and $SO(n)$, and use these to
calculate the equivariant volumes.  Nevertheless, since an interpretation of these desingularizations as moduli spaces is at
least a priori lacking (see however~\cite{thesishans} for related discussions), determining the fixed point data in the hope of applying a direct localization formula is a non-trivial matter.  
In~\cite{incohnonab} the Kirwan desingularization construction was used to develop a residue formula as for intersection
pairings in the intersection cohomology of a singular GIT quotient. While we are not
directly interested in the intersection cohomology of the Uhlenbeck spaces, we can take a similar approach for
calculating the equivariant volumes of the ADHM spaces for symplectic and  special orthogonal gauge groups, to obtain the equivalent formula for these volumes to the one derived in the previous section for $SU(n)$, by considering a degenerate form on the Kirwan desingularization.
\subsubsection{$Sp(n)$}
Following~\cite{brysan}, we can describe the ADHM construction of the Uhlenbeck space for gauge group $Sp(n)$ as a \HK quotient as follows: look at the diagram of linear maps
\begin{equation}\label{tekening}\xymatrix{& V \ar@(dr,ur)[]_{B}\ar@(dl,ul)[]^{A}\ar@/_/[ddl]_{j}\ar[dr]_{\Phi} &  \\ & 
& V^{*}
\\
 W\ar[dr]_{J} & &
 \\
& W^{*}\ar@/_/[uur]_{j^*} & 
}
\end{equation}
where $\Phi$ is a (fixed) \emph{real structure} on $V$, i.e. an isomorphism $V\overset{\Phi}{\rightarrow} V^*$ such that
$\Phi^*=\Phi$, and $J$ is a (fixed) \emph{symplectic structure} on $W$ ($W\underset{\cong}{\overset{J}{\rightarrow}}W^*,\
J^*=-J$).   The space of ADHM data \sympladhm{n}{c} here consists of $\{(A,B,j)\}$, with the extra conditions that $\Phi A, \Phi B\in S^2 V^*$, and the group divided by is $O(V)$, determined by $\Phi$.  We can write the vanishing of the complex moment map as $\Phi [A,B] - j^* J j =0$, and  we can again put
together a monad:

\begin{equation}\label{monadski}
V\otimes \mathcal{O}(-1) 
\xrightarrow{ {{\displaystyle a=}{
\left(\! \! {\renewcommand{\arraystretch}{0.4}
\begin{array}{c} \scriptscriptstyle x_0 A-x_1 \\ \scriptscriptstyle x_0 B-x_2 \\ \scriptscriptstyle x_0 j 
\end{array} }\! \!
\right)
}}}
{\renewcommand{\arraystretch}{0.7}
\begin{array}{c}
V\otimes \mathcal{O} \\ \oplus  \\ V\otimes \mathcal{O}\\ \oplus \\ W\otimes \mathcal{O} 
\end{array}}
\xrightarrow{{\displaystyle b=}
\left(\!  {\renewcommand{\arraystretch}{0.4}
\begin{array}{ccc} \scriptscriptstyle x_2 \Phi -x_0 B^* \Phi  & \scriptscriptstyle - x_1 \Phi + x_0 A^* \Phi & \scriptscriptstyle - x_0  j^* J 
\end{array} } \! 
\right)
 }
 V\otimes \mathcal{O}(1).
\end{equation}
Composing this with $\Phi$ on the right gives a self-dual
monad, whose cohomology is an $Sp(n)$ instanton.

As in the case of $SU(n)$, we are again interested in a torus action given by changing the framing at infinity and
rescaling $\C^2\subset \CP_2$.  The former action is readily lifted to \sympladhm{n}{c}: 
$$t\in T_{Sp(n)}\Rightarrow t. (A,B,j)=(A,B,t j).$$
As for the scaling, things become a bit more cumbersome --- a
lift of the torus action to the space of ADHM data does not seem to be available.  However, we can still proceed as before if we temper our
ambition and only try to lift the action with weight two --- that is to say the action induced by the scaling of
$\C^2\subset \CP_2$ given by
$$(e_1,e_2).(x_0,x_1,x_2)=(x_0, e_1^2x_1, e_2^2x_2).$$
Indeed, if we now introduce the bundle isomorphism 
$$\phi:{\renewcommand{\arraystretch}{0.8}
\begin{array}{c}
V\otimes \mathcal{O} \\ \oplus  \\ V\otimes \mathcal{O}\\ \oplus \\ W\otimes \mathcal{O} 
\end{array}} \longrightarrow {\renewcommand{\arraystretch}{0.8}
\begin{array}{c}
V\otimes \mathcal{O} \\ \oplus  \\ V\otimes \mathcal{O}\\ \oplus \\ W\otimes \mathcal{O} 
\end{array}}: \left(\! {\renewcommand{\arraystretch}{0.9}
\begin{array}{c}
v_1 \\ v_2 \\ w
\end{array}}
 \!\right)\mapsto 
 \left(\!{\renewcommand{\arraystretch}{1.3}
\begin{array}{c}
\frac{e_1}{e_2} v_1 \\ \frac{e_2}{e_1} v_2 \\w
\end{array}}
\! \right),$$
then by using this isomorphism we obtain
$$\im \left(\! \! {\renewcommand{\arraystretch}{0.4}
\begin{array}{c} \scriptstyle x_0 A-e_1^{-2}x_1 \\ \scriptstyle x_0 B-e_2^{-2}x_2 \\ \scriptstyle x_0 j 
\end{array} }\! \! 
\right)
\stackrel{\phi}{\cong} \im \frac{1}{e_1 e_2}\left( \! \!  {\renewcommand{\arraystretch}{0.4}
\begin{array}{c} \scriptstyle x_0 (e_1^2A)- x_1 \\ \scriptstyle x_0 (e_2^2 B) -x_2 \\ \scriptstyle x_0 (e_1e_2 j) 
\end{array} }\! \!
\right)$$
 and   $$\begin{array}{c} \ker \left( {\renewcommand{\arraystretch}{0.4}
\begin{array}{ccc} \scriptstyle e_2^{-2}x_2 \Phi -x_0 B^*\Phi & \scriptstyle -e_1^{-2}x_1\Phi+ x_0 A^*\Phi &
\scriptstyle x_0 j 
\end{array} }
\right) \\
 {\rotatebox[origin=c]{-90}{\ensuremath{\cong}}} \ \phi \\
\ker \frac{1}{e_1e_2} \left({\renewcommand{\arraystretch}{0.4}
\begin{array}{ccc} \scriptstyle x_2\Phi -x_0 (e_2^2B)\Phi &  \scriptstyle -x_1\Phi+x_o(e_2^2B)^*\Phi
& \scriptstyle x_0 (e_1e_2j) 
\end{array} } 
\right).
\end{array}
$$
Hence we can lift the scaling to $$(e_1,e_2).(A,B,j)=(e_1^2A, e_2^2 B, e_1e_2 j).$$
As we are just interested in calculating the equivariant volumes there is no problem in lifting a `higher weight' --- as
the equivariant volumes for the different weights are related by a scaling of the Lie algebra.  Again we can remark that
there is a $U(1)$ in $T$ whose action is given by a constant
global weight (two, in this case):
$$U(1)\hookrightarrow T=T^2\times T_{Sp(n)}:s\mapsto (s,s,1).$$
If we now apply the symplectic cut with respect to this circle action --- using a weight 2 action however on the copy of $\C$ used in the symplectic cut --- to compactify we again get a projective space,
$$\csympladhm{n}{c}=\sympladhm{n}{c}\coprod \mathbb{P}\sympladhm{n}{c},$$ and we can implement the method as before.    Two remarks here: there is an extra factor $\frac{1}{2}$ needed to account for the fact that the moduli space is the \HK quotient by $O(V)$ rather than just $SO(V)$ (the residue formula assumes that the group is connected), and secondly, because we only lifted the `weight 2' action the space of ADHM data the corresponding variables $\epsilon_1, \epsilon_2$ need to be rescaled.  

For the calculations below we will fix isomorphisms $V\cong \C^c$, $W\cong \C^{2n}$, and represent $\Phi$ by the off-diagonal matrix $\left({\renewcommand{\arraystretch}{0.1}
\begin{array}{ccc}  0  &  &    \ \ \scriptstyle 1  \\ & \! {\rotatebox[origin=c]{45}{\ensuremath{\cdots}}} \! \! \!
   \! \!  & \\ \!\! \scriptstyle 1 \! & \! \! \!   \! \! \! & \ \ \ 0
\end{array} }
\right)$, and $J$ by $\left({\renewcommand{\arraystretch}{0.8}
\begin{array}{cc}  0 & \mathbb{I}_n  \\  -\mathbb{I}_n &  0
\end{array} }
\right)$.  Using maximal tori of the form $\text{diag}(s_1,\ldots, s_m,s_m^{-1},\ldots,s_1)$ for even $c=2m$, $\text{diag}(s_1,\ldots,s_m,1,s_m^{-1},\ldots,s_1^{-1})$ for odd $c=2m+1$ for $SO(c)$ and  $\text{diag}(t_1,\ldots,t_n,t_1^{-1},\ldots,t_n^{-1})$ for $Sp(n)$, we can easily write down all the weights involved.  Implementing the calculational method we get
$$\int_{\sympl{n}{c}} e^{\omega+\mu} = \frac{1}{2}\frac{1}{|W|} \res_+^{\sigma_i} \frac{\varpi^2 \ e_T(\mu_{\C})}{e_T(\sympladhm{n}{c})},$$
where for $c=2m$ even
$$|W|=2^{m-1}m!$$
$$\varpi^2=\prod_{i<j} (\sigma_i^2-\sigma_j^2)^2$$
$$e_T(\mu_{\C})=\left(\epsilon_1+\epsilon_2\right)^m  \prod_{i<j} \left(\left(\epsilon_1+\epsilon_2\right)^2-\left(\sigma_i+\sigma_j\right)^2\right)\left(\left(\epsilon_1+\epsilon_2\right)^2-\left(\sigma_i-\sigma_j\right)^2\right)$$
and 
\begin{multline*}e_T(\sympladhm{n}{c})=\prod_{k=1,2}\left(\prod_{i,j}(\epsilon_k +\sigma_i -\sigma_j)\prod_{i\leq j}\left((\epsilon_k)^2-(\sigma_i+\sigma_j)^2\right)\right)\\ \prod_i\prod_l\left( \left(\frac{\epsilon_1+\epsilon_2}{2}+\tau_l\right)^2-\sigma_i^2 \right) \left( \left(\frac{\epsilon_1+\epsilon_2}{2}-\tau_l\right)^2-\sigma_i^2 \right),\end{multline*}
and for $c=2m+1$ odd
$$|W|=2^{m}m!$$
$$\varpi^2=\prod_{i<j} (\sigma_i^2-\sigma_j^2)^2\prod_i \sigma_i^2$$
\begin{multline*}e_T(\mu_{\C})=\left(\epsilon_1+\epsilon_2\right)^{m}\prod_i\left((\epsilon_1+\epsilon_2)^2-\sigma_i^2\right)\\ \prod_{i<j} \left(\left(\epsilon_1+\epsilon_2\right)^2-\left(\sigma_i+\sigma_j\right)^2\right)\left(\left(\epsilon_1+\epsilon_2\right)^2-\left(\sigma_i-\sigma_j\right)^2\right)\end{multline*}
and 
\begin{multline*}e_T\left(\sympladhm{n}{c}\right)= \prod_{k=1,2}\!\!\left(\!\!\epsilon_k\!\prod_i(\epsilon_k^2-\sigma_i^2)\!\prod_{i,j}(\epsilon_k+\sigma_i -\sigma_j)\!\prod_{i\leq j}\!\left((\epsilon_k)^2-(\sigma_1+\sigma_2)^2\right)\!\!\right)\\ \prod_i \prod_l \! \! \left( \! \!\left(\frac{\epsilon_1+\epsilon_2}{2}+\tau_l\right)^2\! \!-\sigma_i^2 \right) \! \! \left(\! \! \left(\frac{\epsilon_1+\epsilon_2}{2}-\tau_l\right)^2\! \!-\sigma_i^2 \right)\! \! \prod_l \! \!\left(\! \!\left(\frac{\epsilon_1+\epsilon_2}{2}\right)^2\! \! - \tau_l^2\right).\end{multline*}

\subsubsection{SO(n)}
For the gauge group $SO(n)$ the ADHM construction goes similarly: again using the diagram (\ref{tekening}), where in this case  $\Phi:V\rightarrow V^{*}$ is a symplectic structure,  and
$J:W\rightarrow W^{*}$ is a real structure,
the space of ADHM data $\sportadhm{c}{n}$ consists of the
triples of linear maps $(A,B,j)$, this time with the extra
condition that $\Phi A,\ \Phi B \in \bigwedge^2 V^*$.  
With this understood, (and the vanishing of the complex moment map again being $\Phi[A,B]-j^*Jj=0$) nominally the same monad (\ref{monadski}) again gives the desired bundle, from which we can directly see that the same lift of the torus action works (again with `double weight' for the 2-torus that scales): $$(e_1,e_2,t).(A,B,j)=(e_1^2 A, e^2_2 B,e_1e_2tj)$$ for $(e_1,e_2)\in T^2$ and $j\in T_{SO(n)}$.

With maximal tori written as before it is again just a matter of identifying the weights for all tori involved on the space \sportadhm{n}{c} and the complex moment map.  We obtain again
$$\int_{\sport{n}{c}} e^{\omega+\mu} = \frac{1}{c!2^c} \res_+^{\sigma_i} \frac{\varpi^2 \ e_T(\mu_{\C})}{e_T(\sportadhm{n}{c})},$$
where
$$\varpi^2=\prod_{i<j} (\sigma_i^2-\sigma_j^2)^2\prod_i (2\sigma_i)^2$$
$$e_T\left(\mu_{\C}\right)= (\epsilon_1+\epsilon_2)^{\frac{n}{2}}\prod_{i<j}\left((\epsilon_1+\epsilon_2)^2-(\sigma_i-\sigma_j)^2\right)  \prod_{i\leq j}\left( (\epsilon_1+\epsilon_2)^2-(\sigma_i+\sigma_j)^2\right), $$
 for even $n=2m$
\begin{multline*}e_T(\sportadhm{n}{c})=\prod_{k=1,2}\left(\prod_{i,j} (\epsilon_k+\sigma_i-\sigma_j)\prod_{i<j} (\epsilon_k^2-(\sigma_i+\sigma_j)^2)\right)\\ \prod_i\prod_l\left( \left(\frac{\epsilon_1+\epsilon_2}{2}+\tau_l\right)^2-\sigma_i^2 \right) \left( \left(\frac{\epsilon_1+\epsilon_2}{2}-\tau_l\right)^2-\sigma_i^2 \right),\end{multline*}
and for odd $n=2m+1$
\begin{multline*}e_T(\sportadhm{n}{c})=\!\!\prod_{k=1,2}\!\!\left(\!\prod_{i,j} (\epsilon_k\!+\!\sigma_i\!-\!\sigma_j)\prod_{i<j} (\epsilon_k^2\!-(\sigma_i+\sigma_j)^2)\!\!\right)\!\!\prod_i\!\! \left(\!\!\left(\frac{\epsilon_1+\epsilon_2}{2}\right)^2\!\!\!-\sigma_i^2\right)\\ \prod_i\prod_l\left( \left(\frac{\epsilon_1+\epsilon_2}{2}+\tau_l\right)^2-\sigma_i^2 \right) \left( \left(\frac{\epsilon_1+\epsilon_2}{2}-\tau_l\right)^2-\sigma_i^2 \right).\end{multline*}
The above formulas exactly correspond to those found in~\cite[\S 5.3]{ABCD}.

\subsection{Comparison}
We will now compare the calculational strategy outlined above with the work of Nekrasov-Shadchin~\cite{ABCD}.

In~\cite[\S4]{naka} a $K$-theoretic interpretation of the
coefficients of Ne\-kra\-sov's partition function
(\ref{nekrapart}) is given as follows.  Let $E$ be a $T$-equivariant coherent sheaf on $\gies{n}{c}$ and look at 
\begin{equation}\label{ktheo}\sum_{i=0}^{2nr} (-1)^i \ch H^i
  (\gies{n}{c},E) \end{equation} where the character $\ch$
denotes the trace of the representation, defined as a Hilbert
series for a representation $M$ with weight decomposition $M=\bigoplus_{\mu}M_{\mu}$  as
$$\ch M=\sum_{\mu} \dim M_{\mu} e^{\mu}$$ where the sum is over all characters.  For the case of the structure sheaf $\mathcal O$,
(\ref{ktheo}) reduces to $\ch H^0(\uhl{c}{n},\mathcal{O})$, and the following $K$-theoretic localization formula is explained~\cite[Prop. 4.1]{naka}:
$$ \sum_i (-1)^i \ch H^i (\gies{n}{c}, E)=\sum_{F\subset \left(\gies{n}{c}\right)^{T}} \frac{ i^*_F E}{\bigwedge_{-1}\nu_F^*},$$
where as usual the sum in the right-hand side is over all fixed points, and $\bigwedge_{-1}$ stands for the alternating sum of exterior powers as an element in the equivariant $K$-group $K^{T}$  (for all details see ~\cite{naka}).  For the case of the structure sheaf $\mathcal O$ this reduces to $$\ch H^0 (\uhl{n}{c}, \mathcal{O})=\sum_{F\subset \left(\gies{n}{c}\right)^{T}} \frac{ 1}{\bigwedge_{-1}\nu_F^*}.$$ 
In the (co)-homological setting, in the mathematics literature on instanton counting~\cite{naka,naka3,braver1} most authors use the equivariant integral $\int 1$ as the coefficients for the generating function, where the integral is defined in practice through a localization theorem $\int 1= \sum_F\frac{1}{e(\nu_F)}$.  In~\cite{naka} the link between this and the $K$-theoretic approach is made by
$$\int_{\gies{n}{c}}1 = \sum \frac{1}{e(\nu F)} = \lim_{\beta\rightarrow 0} \sum \frac{\beta^{2nc}}{\bigwedge_{-1}(\nu_F^*)}=\lim_{\beta\rightarrow 0} \beta^{2nc}\ch H^0.$$
One can think of the limit $\beta\rightarrow 0$ as formally
inducing a multiplication by the inverse Todd class prescribed by a Riemann-Roch theorem; this has the effect of changing the denominators of the localization theorem for $K$-theory to those of (co)homology.  

In~\cite{ABCD}, Nekrasov-Shadchin use the same philosophy to
calculate the partition functions for all classical groups.
They physically interpret the limit $\beta\rightarrow 0$ as
arising from considering the $4$-dimensional theory as the
limit of a $5$-dimensional theory compactified on a circle of
radius $\beta$.  Then they calculate $\ch H^0(\uhl{n}{c})$ by
means of the ADHM construction as follows.  Let $\rho_1,\dots,\rho_.$ be the characters of the $T\times T_G$ action on $\adhm{c}{n}$. Then the $T\times T_G$-equivariant character is
$$\ch H^0 (\adhm{n}{c},\mathcal{O}) = \frac{1}{\prod(1- \rho_i)}.$$
Furthermore, the $T\times T_G$ action on  $\mu_{\C}$ is
homogeneous, say with weights $\nu_1,\dots,\nu_L$ for $L=\dim
G$.  Hence we have
$$\ch H^0 (\mu_{\C}^{-1}(0), \mathcal{O})=\frac{\prod (1-\nu_i)}{\prod (1-\rho_i)}.$$
In order to get the character over $\mu_{C}^{-1}(0)\git G$, we need to get the $G$-invariant part of this.  The projection onto the $G$-invariant part is given by averaging over the whole of $G$,
$$\ch H^0 (\mu_{\C}^{-1}(0)\git G, \mathcal{O}) = \frac{1}{\volu G} \int_G \ch H^0 (\mu_{\C}^{-1}(0), \mathcal{O}),$$ but by the Weyl integral formula this can be reduced 
 to an integral over the maximal torus,
 $$\ch H^0 (\mu_{\C}^{-1}(0)\git G, \mathcal{O}) = \frac{1}{\abs{W_G}\volu T_G} \int_{T_G} \ch H^0 (\mu_{\C}^{-1}(0), \mathcal{O}). $$
 Now factor the torus as $T=\prod U(1)$ and break up the
 integral into circle integrals.  The integrand in each lives
 on $\C^*$,  and one can think of these as contour integrals in
 the plane.  These contour integrals can be calculated by
 Cauchy's theorem, for each integral keeping the other
 variables as generic constants.  Finally, take the limit
 $\beta\rightarrow 0^+$ to reduce everything to cohomology.
 Nekrasov-Shadchin interpret this limit as a contour integral
 of a meromorphic top degree form over the complexified Lie
 algebra.  In either case the actual evaluation happens by means of iterated residues, for each variable choosing a half-space in which to consider the poles.  It is thus that the formula gives exactly the same result as our method outlined above.  It is interesting to remark that the non-compactness manifests itself in different ways in the two methods.  
 
Let us illustrate this by calculating the regularized volume of the simplest \HK quotient, $\C^4\hk \C^*$, where $s\in\C^*$ acts as $$s.(x_1, x_2 , x_3, x_4)= (s x_1, s x_2, s^{-1} x_3, s^{-1} x_4).$$  Furthermore we consider a $T^2$ action on $\C^4$ by $$(t_1,t_2)(x_1,x_2,x_3,x_4)= (t_1x_1, t_2 x_2, t_2 x_3, t_1 x_4).$$  As the complex moment map $\mu_{\C}$ for the $\C^*$ action is $x_1x_3+ x_2 x_4$, the $T^2$ action indeed preserves $\mu^{-1}_{\C}(0)$, and furthermore the moment map for the $T^2$ action on $\C^4$ clearly has a component that is proper and bounded below, hence all the conditions are satisfied.   For the calculation according to Nekrasov-Shadchin
we need to find the Hilbert-series of $\mathcal{O}$ on $\mu_{\C}^{-1}(0)$, which is
$$\frac{1-t_1t_2}{(1-st_1)(1-st_2)(1-t_2 s^{-1})(1-t_1s^{-1})}.$$
In order to calculate the Hilbert series over the quotient
$\mu^{-1}_{\C}(0)\git \C^{*}$ we then have to take the contour integral
\begin{equation*}
\begin{split}
\ch &=\frac{1}{2\pi i}\oint \frac{ds}{s} \frac{1-t_1t_2}{(1-st_1)(1-st_2)(s-\frac{t_2}{s})(1-\frac{t_1}{s})} \\ & = \frac{ 1-t_1t_2}{(1-t_1t_2)(1-t_2^2)(1-\frac{t_1}{t_2})}+\frac{1-t_1t_2}{(1-t_1^2)(1-t_1t_2)(1-\frac{t_2}{t_1})}.
\end{split}
\end{equation*}
Finally set $t_1=e^{-\beta \tau_1},\ t_1=e^{-\beta \tau_2}$ and $s=e^{-\beta \sigma}$ and take $\lim_{\beta\rightarrow 0} \beta^2 \ch$, which gives \begin{equation}\label{nekras}\frac{\tau_1+\tau_2}{(\tau_2+\tau_1)(2\tau_2)(\tau_1-\tau_2)}+\frac{\tau_1+\tau_2}{2\tau_1(\tau_1+\tau_2)(\tau_2-\tau_1)}=\frac{1}{2\tau_1\tau_2}.\end{equation}

On the other hand, in the method we have outlined above, we can
use the symplectic cut on $\C^4$ with respect to the diagonal
$T^1\subset T^2$, which gives $(\C^4)_{\lambda}=\CP_3\coprod
\C^4$.  For clarity we shall not take the calculational
shortcut discussed in section~\ref{strat}, but actually go
through the procedure of making the symplectic cut and taking
the limit as $\lambda\rightarrow\infty$.  There are five fixed points for the $\C^*\times T^2$ action on this space, of which we only need to consider $[1\!:\!0\!:\!0\!:\!0]$ and $[0\!:\!1\!:\!0\!:\!0]$ for the Jeffrey-Kirwan residue formula.  Applying this gives 
\begin{multline*}\jkres^{+}_{\sigma} \frac{e^{\lambda(\sigma+\tau_1)}(-2\sigma+\tau_2-\tau_1)}{(-\sigma-\tau_1)(\tau_2-\tau_1)(\tau_2-\tau_1-2\sigma)(-2\sigma)}\\+\jkres_{\sigma}^{+}\frac{e^{\lambda(\sigma+\tau_2)}(-2\sigma+\tau_1-\tau_2)}{(\tau_1-\tau_2)(-\sigma-\tau_2)(-2\sigma)(\tau_1-\tau_2-2\sigma)}. \end{multline*}
As usual, we are only interested in the terms that survive the $\lambda\rightarrow \infty$ limit on the open cone of $Lie(T^2)$; hence for each of the two above terms we only need to consider the residue at the pole that cancels the exponent, which gives the exact same expression as (\ref{nekras}).

\appendix
\section{Proof of theorem~\ref{equiJK}}\label{appen}
On could essentially try to adapt any proof of
theorem~\ref{equiJK} to the equivariant setting.  The approach we take here is almost completely based on~\cite{locbycuts}, and we restrict ourselves to commenting on how to adjust it.
\begin{proof}
We can break up the proof in two steps: first reduce the case where $G$ is a general compact connected Lie group to the
corresponding statement for a maximal torus $T_{G}$ of $G$, and then prove the theorem just for the case when $G$ is a
torus.  The first step is achieved by the abelianization theorem: 

\begin{theorem} \label{equimartin} Let $M$ as above be equipped with commuting Hamiltonian group actions of compact connected Lie groups $G,H$, let $T_G$ be a maximal torus of $G$, and let $a$ be an element of $H^*_{G\times H} M$.  Then in $H^*_H(*)$,
$$\int_{M\git G} \kappa_G (a) =\frac{1}{\abs{W_G}} \int_{M\git T_G} \kappa_{T_G} (\varpi^2 a).$$
\end{theorem}

This theorem is a corollary of the Jeffrey-Kirwan residue theorem~\cite{nonabloc}, but was also proven directly by Martin~\cite{martin}.
The original theorem of Martin was not formulated in an
equivariant setting (i.e.\ no $H$ was present), but it suffices to
remark here that the proof given in~\cite{martin} is valid without modifications in the equivariant case as well.

Now, to prove theorem~\ref{equiJK} for a torus $G=T$, we will
first outline the strategy of the proof given
in~\cite{locbycuts}, and then comment on how it can be used for
our equivariant purposes.  Begin by choosing a cone $\Lambda$,
as follows.  Look at the weights of the torus $T$ at all the components of $M^{T}$.  Each of these determines a hyperplane in $\mathfrak{t}$.  Choose a connected component of the complement of their union; call this $\Lambda$.  Then look at a polyhedral cone containing the dual cone, say $\Lambda^*\subset\Sigma$.  
The next step is to take the symplectic cut $M_{\Sigma}$ of $M$ with respect to  $\Sigma$ --- having chosen $\Sigma$ in such a way that this is smooth.  For technical purposes we here first alter our moment map for the $T$-action a little by choosing a point $p\in \Sigma$, and replacing $\mu_T$ by $\mu_{\epsilon}=\mu_T-\epsilon p$.   As we inherit a commuting action of $T$ and $H$ on this, we can investigate the fixed-point set $M_{\Sigma}^{T}$. The fixed-point components come in three flavors:  `old ones,' i.e. components $F$ with $\mu(F)\in \Sigma^0$, `new ones,' i.e. components $\tilde{F}$ with $0\neq \mu(\tilde{F})\in \partial(\Sigma)$, and the component at the vertex of $\Sigma$, which we can identify with $(M\git_{\epsilon p} T)$.
Applying the localization theorem 
to $M_{\Sigma}$ gives
\begin{equation}\label{misschien}
\int_{M_{\Sigma}}e^{\tilde{\omega}}\ \eta = \sum_F \int_F \frac{e^{\tilde{\omega}}\ i_F^* \eta}{e_{T} (\nu_F)} +\sum_{\tilde{F}} \int_{\tilde{F}} \frac{e^{\tilde{\omega}} i_{\tilde{F}}^* \eta_{\Sigma}}{e_{T_1} (\nu_G)} + \int_{M\git_{\epsilon p} T}\frac{e^{\tilde{\omega}}i^* \eta}{e_{T_1} (\nu_{M\git T})}. 
\end{equation}
As a function on $\mathfrak{t}\oplus \mathfrak{h}$ this is `holomorphic' (i.e. the analytic continuation is holomorphic on the complexification $\mathfrak{t}_{\C}\oplus \mathfrak{h}_{\C}$), and for such functions it is proven by Jeffrey and Kogan~\cite[Lemma 3.3]{locbycuts} that the residue with respect to $\Lambda$ and $-\Lambda$ give the same results.  This leads to the equality
\begin{multline}\label{eindelijk} \JKres^{\Lambda} \left( \sum_F \int_F \frac{e^{\tilde{\omega}}\ i_F^* \eta}{e_{T} (\nu_F)}+\sum_{\tilde{F}}  \int_{\tilde{F}} \frac{e^{\tilde{\omega}}\ i_{\tilde{F}}^* \eta_{\Sigma}}{e_{T_1} (\nu_G)} + \int_{M\git_{\epsilon p} T}\frac{e^{\tilde{\omega}}\ i^* \eta}{e_{T_1} (\nu_{M\git T})} \right)[dx] \\  = \JKres^{-\Lambda} \left( \sum_F \int_F \frac{e^{\tilde{\omega}}\ i_F^* \eta}{e_{T} (\nu_F)} +\sum_{\tilde{F}} \int_{\tilde{F}} \frac{e^{\tilde{\omega}}\ i_{\tilde{F}}^* \eta_{\Sigma}}{e_{T_1} (\nu_G)} + \int_{M\git_{\epsilon p} T}\frac{e^{\tilde{\omega}}\ i^* \eta}{e_{T_1} (\nu_{M\git T})} \right)[dx].  \end{multline}
The modification of the moment map from $\mu_{T}$ to
$\mu_{\epsilon}$ is done exactly to ensure that the residue operation is valid for all the terms appearing above.  Next it is shown in~\cite[\S 5.2]{locbycuts} that in fact, for small enough $\epsilon > 0$ and with the choices of $\Lambda$, $\Sigma$ and $\mu_{\epsilon}$ as above, $4$ terms in this equality are zero, leading to
$$\JKres^{-\Lambda}   \int_{M\git_{\epsilon p} T}\frac{e^{\tilde{\omega}}\  i^* \eta}{e_{T_1} (\nu_{M\git T})} [dx]= \JKres^{\Lambda} \sum_F \int_F \frac{e^{\tilde{\omega}}\ i_F^* \eta}{e_{T} (\nu_F)}[dx]. $$
Finally, take the limit as $\epsilon \to 0$.  By identifying
the weights of $T$ on $\nu_{M\git T}$ (which is a bundle of rank $\dim \mathfrak{t}$), 
one obtains
\begin{equation}\label{pfst}\lim_{\epsilon\rightarrow 0^+} \JKres^{-\Lambda} \int_{M\git T} \frac{e^{\tilde{\omega}}\ i^* \eta }{e_{T_1}(\nu_{M\git T})}[dx]= \int_{M\git T} e^{\omega}\ \kappa{\eta},\end{equation}
and finally this gives 
$$\int_{M\git T} \kappa(\eta)e^{\omega} = \sum_{F\subset M^{T}} \JKres^{\Lambda}\int_{F}\frac{e^{\tilde{\omega}} i^*_F \eta}{e_T(\nu_F)}[dx].$$
This is how theorem~\ref{lisafrances} is proved in the non-equivariant case in~\cite{locbycuts}.
Now in order to make this valid in the presence of the extra
group action of $H$, we remark that because the actions of $T$
and $H$ on $M$ commute, the action of $H$ descends to
$M_{\Sigma}$.  Now we can use our approximation of the Borel model by smooth finite-dimensional spaces $M\times E_iH/H$, determined by finite-dimensional approximations $E_iH\rightarrow B_iH$ to the classifying space $EH\rightarrow BH$.  These
spaces are not symplectic, but they still are Poisson manifolds, with the fibers of $M\times_{H}E_iH\rightarrow B_iH$
as the symplectic leaves.  Furthermore, as the moment map $\mu_{G}:M\rightarrow \mathfrak{g}^{*}$-action is invariant
under $H$, we get moment maps $M\times_{G}E_iH\rightarrow
\mathfrak{g}^{*}$ in the sense of Poisson geometry; see
e.g.~\cite{reduc}.
Also, as the proof of the localization theorem used in (\ref{misschien}) essentially only uses the functoriality of the integral as a push-forward (see~\cite{atibott}), we can apply it to other (proper)
push-forwards as well.  In particular, consider $$
\left(M_{\Sigma} \times E_iH\right) / H \rightarrow B_iH.$$
Denote the corresponding push-forward by $\int_i$, so that $$
\int_i:H^*_{T_G} \left(\left((M_{\Sigma} \times E_iH\right) / H\right) \rightarrow H^*_{T_G}(B_iH)= H^*_{T_G}(*)\otimes
H^*(B_iH)$$
and $$\int_i : H^*_{T_G}\left( \left(F \times E_iH\right)/H\right) \rightarrow H^*_{T_G}(*)\otimes H^*(B_iH).$$
The equivalent of the localization formula is then
$$\int_{i} \alpha = \sum_{\frac{F\times E_iH}{H}\subset \frac{M^T\times E_iH}{H}} \int_i \frac{i^*_F\alpha}{e_T(\nu_{\frac{F\times E_iH}{H}})}.$$
We can now apply the residue operation $\JKres^{\Lambda}$ to both sides (or rather the obvious extension of the residue
from $k(\mathfrak{t})$ to the various $k(\mathfrak{t})\otimes
H^*(B_iH)$), and we observe that for the same reasons
as in~\cite{locbycuts}, the equivalences of (\ref{eindelijk}) and the cancellations of the four terms are still valid ---
this time interpreted as happening in $H^*(B_iH)$.  The same is the case for the equivalent statement of (\ref{pfst}): 
essentially the only thing we need to remark for this is that the weights for the action of $T$ occurring on the normal
bundles of $M\git T$ in $M_{\Sigma}$ are the same as the weights for the normal bundle of  $\frac{M\git T \times E_iH}{
H}$ in $\frac{M_{\Sigma}\times E_iH}{H}$. Furthermore, as all the spaces are finite-dimensional, the inverses of the
Euler classes exist in the usual way.
So we have 
$$\int_{i,M} \kappa (\eta) = -\JKres^{\Lambda}\sum_{F\subset M^T}\int_i \frac{\eta}{e_T(\nu_{\frac{F\times E_iH}{H}})} [dx].$$ 
Now we want to take the limit $i\rightarrow \infty$ to obtain the desired $H$-equivariant result.  We need to observe
in which ring we allow all of the manipulations to take place.  Begin by remarking that in the proof of the abelian localization theorem~\cite{atibott}, which was the basis of all of our localization results, one had to kill some of the torsion in the equivariant cohomology ring, but
tensoring with the full fraction field
$k(\mathfrak{t})$ is a bit excessive.  In particular, we can do with just inverting all products of linear factors
$\prod_i \beta_i$, with $\beta_i\in \mathfrak{t}^*$.  Doing this has the advantage that we still have a grading  on the
ring we obtain (say $\widehat{H^*_T(*)}$).  This gives a bigrading on $\widehat{H^*_{T\times
H}(F)}=\widehat{H^*_T(*)}\otimes H^*_H{F}$ for all $F\subset M^T$.  Now complete these rings by means of the filtration
$$
\dots \subset R_{k+1}\subset R_{k} \subset R_{k-1}\subset \dots$$
where $R_k$ consists of finite sums of elements of bidegree $(i,j)$ such that $ i + 2j \geq k$.  This completion has
all the desired properties (though it is of course not unique as such): it contains the formal cohomology classes
$e^{\tilde{\omega}}$, and it allows us to invert the
equivariant Euler classes --- which because of the splitting principle
we can write without loss of generality as $e_{T\times H}(\nu_F)=\prod_i (\beta_i + c_1 \nu_{i,F})$ with
$\beta_i\in\mathfrak{t}$, $c_1 \nu_{i,F}\in H^2_H{F}$ --- as
$$
\frac{1}{e_{T\times H}(\nu_F)} =
\prod_i\frac{1}{\beta_i}\sum_{j=0}^{\infty}\left(-\frac{c_1(\nu_{i,F})}{\beta_i}\right)^j.$$
The residue $\JKres^{\Lambda}$ of an integral of such a class
will lie in the usual completion of $H^*_H(*)$ by degree.
This completes the proof of theorem~\ref{equiJK}.
\end{proof}

\def\cprime{$'$}


\begin{thebibliography}{DHKM02}

\bibitem[AB84]{atibott}
M.~F. Atiyah and R.~Bott.
\newblock The moment map and equivariant cohomology.
\newblock {\em Topology}, 23(1):1--28, 1984.

\bibitem[AHDM78]{ADHM}
M.~F. Atiyah, N.~J. Hitchin, V.~G. Drinfel{\cprime}d, and Yu.~I. Manin.
\newblock Construction of instantons.
\newblock {\em Phys. Lett. A}, 65(3):185--187, 1978.

\bibitem[BE06]{braver2}
A.~Braverman and P.~Etingof.
\newblock Instanton counting via affine {L}ie algebras. {II}. {F}rom
  {W}hittaker vectors to the {S}eiberg-{W}itten prepotential.
\newblock In {\em Studies in Lie theory}, volume 243 of {\em Progr. Math.},
  pages 61--78. Birkh\"auser Boston, Boston, MA, 2006.

\bibitem[Bra04]{braver1}
Alexander Braverman.
\newblock Instanton counting via affine {L}ie algebras. {I}. {E}quivariant
  {$J$}-functions of (affine) flag manifolds and {W}hittaker vectors.
\newblock In {\em Algebraic structures and moduli spaces}, volume~38 of {\em
  CRM Proc. Lecture Notes}, pages 113--132. Amer. Math. Soc., Providence, RI,
  2004.

\bibitem[BS00]{brysan}
Jim Bryan and Marc Sanders.
\newblock Instantons on {$S\sp 4$} and {$\overline{{\bf C}{\rm P}}{}\sp 2$},
  rank stabilization, and {B}ott periodicity.
\newblock {\em Topology}, 39(2):331--352, 2000.

\bibitem[BV82]{bv}
Nicole Berline and Mich{\`e}le Vergne.
\newblock Classes caract\'eristiques \'equivariantes. {F}ormule de localisation
  en cohomologie \'equivariante.
\newblock {\em C. R. Acad. Sci. Paris S\'er. I Math.}, 295(9):539--541, 1982.

\bibitem[BV99]{verbrion}
Michel Brion and Mich{\`e}le Vergne.
\newblock Arrangement of hyperplanes. {I}. {R}ational functions and
  {J}effrey-{K}irwan residue.
\newblock {\em Ann. Sci. \'Ecole Norm. Sup. (4)}, 32(5):715--741, 1999.

\bibitem[CB01]{craw-boev}
William Crawley-Boevey.
\newblock Geometry of the moment map for representations of quivers.
\newblock {\em Compositio Math.}, 126(3):257--293, 2001.

\bibitem[DH82]{DH2}
J.~J. Duistermaat and G.~J. Heckman.
\newblock On the variation in the cohomology of the symplectic form of the
  reduced phase space.
\newblock {\em Invent. Math.}, 69(2):259--268, 1982.

\bibitem[DH98]{dolgagitflips}
Igor~V. Dolgachev and Yi~Hu.
\newblock Variation of geometric invariant theory quotients.
\newblock {\em Inst. Hautes \'Etudes Sci. Publ. Math.}, (87):5--56, 1998.
\newblock With an appendix by Nicolas Ressayre.

\bibitem[DHKM02]{calc-inst}
Nick Dorey, Timothy~J. Hollowood, Valentin~V. Khoze, and Michael~P. Mattis.
\newblock The calculus of many instantons.
\newblock {\em Phys. Rep.}, 371(4-5):231--459, 2002.

\bibitem[Don84]{instgit}
S.~K. Donaldson.
\newblock Instantons and geometric invariant theory.
\newblock {\em Comm. Math. Phys.}, 93(4):453--460, 1984.

\bibitem[Fre05]{thesishans}
Hans-Georg Freiermuth.
\newblock {\em {On partial compactifications of the space of framed vector
  bundles on the projective plane}}.
\newblock PhD thesis, Columbia University, 2005.

\bibitem[GK96]{guilkalk}
Victor Guillemin and Jaap Kalkman.
\newblock The {J}effrey-{K}irwan localization theorem and residue operations in
  equivariant cohomology.
\newblock {\em J. Reine Angew. Math.}, 470:123--142, 1996.

\bibitem[GLS96]{symfib}
Victor Guillemin, Eugene Lerman, and Shlomo Sternberg.
\newblock {\em Symplectic fibrations and multiplicity diagrams}.
\newblock Cambridge University Press, Cambridge, 1996.

\bibitem[JK95a]{jefkirriem}
Lisa~C. Jeffrey and Frances~C. Kirwan.
\newblock Intersection pairings in moduli spaces of holomorphic bundles on a
  {R}iemann surface.
\newblock {\em Electron. Res. Announc. Amer. Math. Soc.}, 1(2):57--71
  (electronic), 1995.

\bibitem[JK95b]{nonabloc}
Lisa~C. Jeffrey and Frances~C. Kirwan.
\newblock Localization for nonabelian group actions.
\newblock {\em Topology}, 34(2):291--327, 1995.

\bibitem[JK97]{nonabloc2}
Lisa~C. Jeffrey and Frances~C. Kirwan.
\newblock Localization and the quantization conjecture.
\newblock {\em Topology}, 36(3):647--693, 1997.

\bibitem[JK98]{arbitrary}
Lisa~C. Jeffrey and Frances~C. Kirwan.
\newblock Intersection theory on moduli spaces of holomorphic bundles of
  arbitrary rank on a {R}iemann surface.
\newblock {\em Ann. of Math. (2)}, 148(1):109--196, 1998.

\bibitem[JK05]{locbycuts}
Lisa Jeffrey and Mikhail Kogan.
\newblock Localization theorems by symplectic cuts.
\newblock In {\em The breadth of symplectic and Poisson geometry}, volume 232
  of {\em Progr. Math.}, pages 303--326. Birkh\"auser Boston, Boston, MA, 2005.

\bibitem[JKKW03]{incohnonab}
Lisa~C. Jeffrey, Young-Hoon Kiem, Frances Kirwan, and Jonathan Woolf.
\newblock Cohomology pairings on singular quotients in geometric invariant
  theory.
\newblock {\em Transform. Groups}, 8(3):217--259, 2003.

\bibitem[Kir84]{thesisfrances}
Frances~Clare Kirwan.
\newblock {\em Cohomology of quotients in symplectic and algebraic geometry},
  volume~31 of {\em Mathematical Notes}.
\newblock Princeton University Press, Princeton, NJ, 1984.

\bibitem[Kir85]{francesdesing}
Frances~Clare Kirwan.
\newblock Partial desingularisations of quotients of nonsingular varieties and
  their {B}etti numbers.
\newblock {\em Ann. of Math. (2)}, 122(1):41--85, 1985.

\bibitem[Ler95]{cuts}
Eugene Lerman.
\newblock Symplectic cuts.
\newblock {\em Math. Res. Lett.}, 2(3):247--258, 1995.

\bibitem[Mar]{martin}
Shaun~K. Martin.
\newblock Symplectic quotients by a nonabelian group and by its maximal torus.
\newblock To appear in \emph{Annals of Math.}, \texttt{arXiv:math.SG/0001001}.

\bibitem[Mei06]{meinrequi}
Eckhard Meinrenken.
\newblock {Equivariant cohomology and the Cartan model}.
\newblock In Jean-Pierre Fran{\c{c}}oise, Gregory~L. Naber, and Tsou~Sheung
  Tsun, editors, {\em Encyclopedia of mathematical physics.} Academic
  Press/Elsevier Science, Oxford, 2006.

\bibitem[MNS00]{higgs-int}
Gregory Moore, Nikita Nekrasov, and Samson Shatashvili.
\newblock Integrating over {H}iggs branches.
\newblock {\em Comm. Math. Phys.}, 209(1):97--121, 2000.

\bibitem[Nak99]{nakabook}
Hiraku Nakajima.
\newblock {\em Lectures on {H}ilbert schemes of points on surfaces}, volume~18
  of {\em University Lecture Series}.
\newblock American Mathematical Society, Providence, RI, 1999.

\bibitem[Nek03]{nekra}
Nikita~A. Nekrasov.
\newblock Seiberg-{W}itten prepotential from instanton counting.
\newblock {\em Adv. Theor. Math. Phys.}, 7(5):831--864, 2003.

\bibitem[NO06]{rndprt}
Nikita~A. Nekrasov and Andrei Okounkov.
\newblock Seiberg-{W}itten theory and random partitions.
\newblock In {\em The unity of mathematics}, volume 244 of {\em Progr. Math.},
  pages 525--596. Birkh\"auser Boston, Boston, MA, 2006.

\bibitem[NS98]{noncominst}
Nikita Nekrasov and Albert Schwarz.
\newblock Instantons on noncommutative {$\mathbf R\sp 4$}, and {$(2,0)$}
  superconformal six-dimensional theory.
\newblock {\em Comm. Math. Phys.}, 198(3):689--703, 1998.

\bibitem[NS04]{ABCD}
Nikita Nekrasov and Sergey Shadchin.
\newblock A{BCD} of instantons.
\newblock {\em Comm. Math. Phys.}, 252(1-3):359--391, 2004.

\bibitem[NY04]{naka2}
Hiraku Nakajima and K{\=o}ta Yoshioka.
\newblock Lectures on instanton counting.
\newblock In {\em Algebraic structures and moduli spaces}, volume~38 of {\em
  CRM Proc. Lecture Notes}, pages 31--101. Amer. Math. Soc., Providence, RI,
  2004.

\bibitem[NY05a]{naka}
Hiraku Nakajima and K{\=o}ta Yoshioka.
\newblock {Instanton counting on blowup, I. 4-Dimensional pure gauge theory}.
\newblock {\em Invent. Math.}, 162(2):313--355, 2005.

\bibitem[NY05b]{naka3}
Hiraku Nakajima and K{\=o}ta Yoshioka.
\newblock Instanton counting on blowup. {II}. {$K$}-theoretic partition
  function.
\newblock {\em Transform. Groups}, 10(3-4):489--519, 2005.

\bibitem[OR04]{reduc}
Juan-Pablo Ortega and Tudor~S. Ratiu.
\newblock {\em Momentum maps and {H}amiltonian reduction}, volume 222 of {\em
  Progress in Mathematics}.
\newblock Birkh\"auser Boston Inc., Boston, MA, 2004.

\bibitem[OSS80]{okonek}
Christian Okonek, Michael Schneider, and Heinz Spindler.
\newblock {\em Vector bundles on complex projective spaces}, volume~3 of {\em
  Progress in Mathematics}.
\newblock Birkh\"auser Boston, Mass., 1980.

\bibitem[Par00]{paradan}
Paul-Emile Paradan.
\newblock The moment map and equivariant cohomology with generalized
  coefficients.
\newblock {\em Topology}, 39(2):401--444, 2000.

\bibitem[PW94]{prato-dh}
Elisa Prato and Siye Wu.
\newblock Duistermaat-{H}eckman measures in a non-compact setting.
\newblock {\em Compositio Math.}, 94(2):113--128, 1994.

\bibitem[SW94]{SW}
N.~Seiberg and E.~Witten.
\newblock Electric-magnetic duality, monopole condensation, and confinement in
  {$N=2$} supersymmetric {Y}ang-{M}ills theory.
\newblock {\em Nuclear Phys. B}, 426(1):19--52, 1994.

\bibitem[Tha96]{gitflips}
Michael Thaddeus.
\newblock Geometric invariant theory and flips.
\newblock {\em J. Amer. Math. Soc.}, 9(3):691--723, 1996.

\bibitem[Ver96]{notevergne}
Mich{\`e}le Vergne.
\newblock A note on the {J}effrey-{K}irwan-{W}itten localisation formula.
\newblock {\em Topology}, 35(1):243--266, 1996.

\end{thebibliography}
\end{document}